\documentclass[a4paper,twoside,12pt]{article}

\usepackage[utf8]{inputenc}
\usepackage{amsmath,amsfonts,amssymb}
\usepackage{vmargin,graphicx,theorem}
\usepackage[english]{babel}
\usepackage{enumerate}
\usepackage{color}
\usepackage{pstricks}
\usepackage{pst-fill,pst-grad,pst-plot,pst-eucl,pstricks-add,pst-node}

\setpapersize[portrait]{A4}
\setmarginsrb{2.5cm}{1.5cm}{2.5cm}{1cm}{0cm}{0.1cm}{0.5cm}{2cm}

\newcommand{\disp}{\displaystyle}

\newcommand{\dN}{\ensuremath{\mathbb{N}}}

\newcommand{\dR}{\ensuremath{\mathbb{R}}}

\newtheorem{ethm}{Theorem}[section]

\newtheorem{eprop}[ethm]{Proposition}

\newtheorem{elem}[ethm]{Lemma}

\newtheorem{edefi}[ethm]{Definition}

\newtheorem{erem}[ethm]{Remark}

\newcommand{\proofend}{~$\rhd$}
\newcommand{\proofbegin}{~$\lhd$}

\newenvironment{eproof}
               {\noindent {{\it \textbf{Proof}}}\\\proofbegin~}
               {\proofend\\}

\newcommand{\ESSUP}[1]{\ensuremath{{\text{ess sup}_{#1}}}}
\newcommand{\PAR}[1]{\ensuremath{{\left(#1\right)}}} 
\newcommand{\BRA}[1]{\ensuremath{{\left\{#1\right\}}}} 
\newcommand{\NRM}[1]{\ensuremath{{\left\Vert #1\right\Vert}}} 

\renewcommand{\geq}{\geqslant}

\newcommand{\ent}{\ensuremath{\text{Ent}}}

\def\disp{\displaystyle}

\newcommand{\A}{\ensuremath{\mathcal A}}
\newcommand{\B}{\ensuremath{\mathcal B}}

\newcommand{\N}{\ensuremath{\mathbb{N}}}

\newcommand{\beq}{\begin{equation}}\newcommand{\eeq}{\end{equation}}
\newcommand{\su}{\mathrm{u}}
\newcommand{\sv}{\mathrm{v}}

\newtheorem{theorem}{Theorem}
\newtheorem{proposition}[theorem]{Proposition}
\newtheorem{lemma}[theorem]{Lemma}

\newtheorem{remark}[theorem]{Remark}

\begin{document}

\title{Solution of a class of  reaction-diffusion systems
 via logarithmic Sobolev inequality}
\author{ Pierre Foug\`eres\thanks{Universit\'e Toulouse 3, IMT}, Ivan Gentil\thanks{Universit\'e Lyon 1, ICJ}\, and Boguslaw Zegarlinski\thanks{Imperial College, London} }

\date{\today}

\maketitle

\abstract{We study global existence, uniqueness and positivity of  
weak solutions 
of a class  of reaction-diffusion systems  of chemical kinetics type, 
under the assumptions of logarithmic Sobolev inequality
and appropriate exponential integrability of the initial data.}

\bigskip

\bigskip

\section{Introduction}

%
%

A mixture one gets after esterification of one mole of ethyl alcohol by one mole of ethanoic acid
contains products (ethyl acetate and water), but also reactants. 
This is an example of a double 
displacement reaction 
\begin{equation}
\label{eq-example}
CH_3CH_2OH + CH_3COOH \rightleftharpoons CH_3COOCH_2CH_3 + H_2O
\end{equation}
(see \cite{mahe-fraissard}).

\noindent
We consider here chemical reactions 
between  $q\geq2$ species $\A_i$, $i=1, \dots, q$, as follows 
\begin{equation*}
\sum_{i=1}^q \alpha_i \A_i \rightleftharpoons \sum_{i=1}^q \beta_i\A_i,
\end{equation*}
where $\alpha_i,\beta_i\in\dN$. 
We assume that for any $1\leq i\leq q$, $\alpha_i-\beta_i\neq 0$ 
which corresponds to the case of a reaction without a catalyst. 

\medskip

If $\vec{u}=(u_1,\cdots,u_q)$ denotes 
the concentration of the species $\A_i$ then the law of action 
mass proposed by Waage and Guldberg in 1864 
{(see again \cite{mahe-fraissard})}
implies that the concentrations are solutions of the system, for all $i\in\{1,\cdots,q\}$,
\begin{equation*}
   \quad\frac{d}{dt}u_i=
   (\beta_i-\alpha_i)\PAR{k\prod_{j=1}^q u_j^{\alpha_j}-l\prod_{j=1}^q u_j^{\beta_j}},
\end{equation*}
where $k,l>0$ are the rate constants of the two reactions. 

\bigskip 

When considering  substances distributed in space, the concentrations 
change 
not only under the influence of the chemical reactions 
but also due to the diffusion of the species over the space, one gets 
the following kinetic model for a chemical reaction-diffusion equation 
$$
\partial_t u_i= L_i u_i +(\beta_i-\alpha_i)\PAR{k\prod_{j=1}^q u_j^{\alpha_j} -l\prod_{j=1}^q u_j^{\beta_j}},
$$
where for all $i=1, \dots, q$,
$L_i$ is an operator which modelizes how the substance diffuses.

We will assume that $L_i=C_i L$ for some $C_i\geq0$ and
some reference operator $L$. 
Moreover, by a change of variables, 
one can assume that there exist (a posteriori two) constants $\lambda_i>0$ 
such that the system of reaction-diffusion is given by 
\begin{equation}
\label{eq-gene}
\partial_t u_i= C_i L u_i +\lambda_i(\beta_i-\alpha_i)\PAR{\prod_{j=1}^q u_j^{\alpha_j} -\prod_{j=1}^q u_j^{\beta_j}},
\end{equation}
where $\vec{u}(t,x)=(u_1(t,x),\cdots,u_q(t,x))$ with $t\geq0$ and $x$ belongs to
the underlying space. 

The  {\it two-by-two} system, one of the simplest non trivial example,
describes the chemical reaction 
$$
\A_1+\A_2 \rightleftharpoons \B_1+\B_2,
$$
and the system of equations can by formulated as follow 
\begin{equation}
\label{eq-tbt}
\left\{
\begin{array}{l}
\partial_t u_1 = C_1L u_1 - \lambda\left(u_1u_{2} - v_1v_{2}\right) \\
\partial_t u_2 = C_2L u_2 -  \lambda\left(u_1u_{2} - v_1v_{2}\right) \\
\partial_t v_1 = C_3L v_1 +\tilde{\lambda}\left(u_1u_{2} - v_1v_{2}\right)\\
\partial_t v_2 = C_4L v_2 + \tilde{\lambda}\left(u_1u_{2} - v_1v_{2}\right)\\
\end{array}
\right.
\end{equation}
where $\lambda, \tilde{\lambda}>0$ and  $u_i$ denotes the concentration of 
the specie $\A_i$ and $v_i$ the concentration of the specie $\B_i$ for $i=1,2$.
To make things even simpler, 
we will assume later that $\lambda=\tilde{\lambda}$.

\bigskip
More general reaction-diffusion systems, of the following form
\begin{equation}
\label{eq-reactiondiffusion}
\left\{
\begin{array}{l}
 \partial_t \vec{u} = \mathfrak{C} \Delta_x \vec{u}+ F(t,x,\vec{u}), \quad 
t >0, \quad x \in \Omega\\
\vec{u}(0)= \vec{u}_0,
\end{array}
\right. 
\end{equation}
with prescribed boundary conditions, were intensively studied in the past.
Here, $\Omega$ is a (possibly unbounded sufficiently smooth) domain of $\mathbb{R}^n$,
$\vec{u}$ takes values in $\dR^q$,
$\mathfrak{C}$ is a usually diagonal $q \times q$ matrix
which can be degenerate, and $F(t,x,\cdot)$ is a vector field on 
$\dR^q$. 

Depending on specific choices for $\mathfrak{C}$ and $F(t,x,\cdot)$,
such systems can present various behaviours with respect to global
existence and asymptotic behaviour of the solution. 
Paragraph 15.4 in \cite{taylorIII} is a nice introduction with a lot
of classical references. 

In the above setting,
local existence follows from general textbooks on parabolic type
partial differential equations (see \cite{friedman}, \cite{ladyzenskaja},
or for fully general boundary value problems \cite{amann-local}).

Global existence question (or how to prevent blow up) gave rise to
extensive efforts and to different methods adapted to specific cases
(see \cite{amann-global}, especially remark 5.4. a), \cite{rothe}, \cite{pierre} and references therein).
Most of these methods consist in deducing $\mathbb{L}^\infty$
bounds on the maximal solution from bounds in weaker norms.

The survey \cite{pierre} provides a lot of references, positive and
negative results, together with a description of open problems.
Its first observation is that, for numerous reaction-diffusion
systems of interest in applications, the nonlinearity satisfies
two general conditions which ensure respectively positivity and
a control of the mass (i.e. the $\mathbb{L}^1$ norm) of a solution.
M. Pierre investigates how these $\mathbb{L}^1$ estimates 
(as well as $\mathbb{L}^1$ bounds on the nonlinearity)
help to provide global existence.

Further works provide results on asymptotic behaviour.
Spectral gap, logarithmic Sobolev
inequality and entropy methods are often used to quantify 
exponential convergence 
of the solution
of an equation to equilibrium, 
and in the context of reaction-diffusion equations (mostly of type \eqref{eq-gene}) they were used to study
the convergence (to constant steady states) in~\cite{df06,df08,desvillettes-seul,gentil-zeg}.  
Geometric characteristics and approximations of global and exponential attractors of general reaction-diffusion systems 
may be found in \cite{zelik03,zelik04,zelik07} (and references therein) in terms 
of precise estimates of their Kolmogorov $\varepsilon$-entropy. In these
papers, $\mathfrak{C}$ is of positive symmetric part and
the nonlinearity must satisfy some moderate growth bound involving the dimension $n$
to ensure global existence. Other cross-diffusion systems are studied by entropy 
methods in \cite{carrillo-hittmeir-jungel}.

\bigskip
One way or another, local or global existence results in the above setting
rely on regularity theory for the heat semigroup, the maximum principle,
and Sobolev inequality through one of its consequences, Gagliardo-Nirenberg
inequalities or ultracontractivity of the semigroup 
(as well as Moser estimates).
(Note nevertheless that an approach based on a nonlinear Trotter product
formula is proposed in \cite{taylorIII}, but seems to impose 
some kind of uniform continuity of the semigroup).

\bigskip 
The aim of this article is to prove global existence of a non-negative 
solution of the reaction-diffusion system~\eqref{eq-gene} with unbounded 
initial data in a setting where Sobolev inequality possibly
does not hold, as e.g. in infinite dimensions or when underlying measure does not satisfy polynomial growth condition. We restrict ourselves to some nonlinearities for 
which in a finite dimensional 
setting, $\mathbb{L}^\infty$ bounds of the solution (and so
global existence) come for free,  \cite{pierre}. Nevertheless,
Sobolev inequality has to be replaced by the weaker logarithmic
Sobolev inequality (or other coercive inequalities which 
survive the infinite dimensional limit; see \cite{bobkov-zegarlinski},
\cite{barthe-cattiaux-roberto}, \cite{roberto-zegarlinski}, 
\cite{bobkov-zegarlinski-slow-tails}).

The celebrated paper \cite{gross75} of L.Gross established
equivalence of logarithmic Sobolev inequality and
hypercontractivity of the semigroup. No compactness embeddings
hold in this context.

For a wide variety of strongly mixing
Markov semigroups, logarithmic Sobolev inequality holds for the corresponding Dirichlet form of the generator.
For diffusion semigroups on Riemannian manifolds, logarithmic 
Sobolev inequality follows from positive bound from below of the Ricci 
curvature (of the generator $L$), 
i.e. the so called Bakry-Emery, 
 $\Gamma_2$ or $\text{CD}(\rho,\infty)$
criterion (see \cite{bakry-emery}, 
\cite{bakry-gentil-ledoux}, \cite{wang}). 
For logarithmic Sobolev inequality and discrete state space Markov chains, 
see \cite{diaconis-saloff}, \cite{miclo}, \cite{chen}.
In infinite dimensional spaces,   
logarithmic Sobolev inequality
for spin systems has been extensively studied 
(see \cite{zegarlinski90}, \cite{stroock-zegarlinski},
\cite{zegarlinski96}, \cite{bodineau-helffer}, \cite{yoshida00}, 
\cite{ledoux-spins}, \cite{guionnet-zegarlinski},  \cite{otto-reznikoff};
and in subelliptic setting, \cite{lugiewicz-zegarlinski},
\cite{inglis-papageorgiou}, \cite{hebisch-zegarlinski}).
%
In the present paper, logarithmic Sobolev
inequality plays  a key role to study existence results
in an infinite dimensional setting, by a step-by-step approximation 
approach.
 

\medskip

The paper is organized as follows. In the next section we describe the framework
and the main result of the paper: in the two-by-two case, assuming 1)~that 
$C_1=C_3$ and $C_2=C_4$, 2)~that the linear {\it diffusion} term satisfies
logarithmic Sobolev inequality and 3)~that the initial datum $\vec{f}$
is nonnegative and satisfies some exponential integrability properties (made more precise
later),  
then there exists a unique  weak solution of the system of 
reaction-diffusion equation~\eqref{eq-tbt} which is moreover nonnegative. 
Section~\ref{section:iterative_approach} presents the iterative procedure
we follow to approximate weak solutions of our reaction-diffusion type problem.
This is based on some cornestone linear problem which is stated there. 
The two following sections are devoted to the details of the proof : section
\ref{section:proof_main} to the convergence of the iterative procedure to
the unique nonnegative weak solution of the nonlinear  Cauchy problem, 
whereas section 
\ref{section:proof_cornerstone} focuses on the cornerstone
linear problem.

In Section~\ref{section:extension} we extend our result 
to the general case of  
system~\eqref{eq-gene}, and present how operators $C_i L$ can be modified. 
(To give a comprehensive proof we focus in the rest of the paper
on the two-by-two case
which already contains non trivial difficulty).

We recall or detail
tools used in the proof in three appendices: the entropic inequality, 
basics  on Orlicz spaces, and finally some further topics on Markov
semigroups and Orlicz spaces.

\section{Framework and main result}
\label{section:main_result}

\subsection*{An abstract Reaction-Diffusion equation}
\label{section:abstract_problem_assumptions}
In the following we will consider an underlying Polish space  
$\mathbb{M}$ equipped with a probability measure $\mu$. 
Let $L$ be a (linear) densely defined selfadjoint Markov operator on $\mathbb{L}^2 (\mu) \equiv \mathbb{L}^2 (\mathbb{M},\mu) $, that is the infinitesimal generator 
of a $C_0$ Markov semigroup $(P_t)_{t \geq 0}$ symmetric with
respect to $\mu$.  It is well known that under these assumptions
there exist a kernel $p_t(x,dy)$ on
$(\mathbb{M}, \mathcal{B}_{\mathbb{M}})$, that is 
a measurable family of probability measures  such that, for any $t \geq
0$, {any $f \in \mathbb{L}^1(\mu)$}, and for $\mu$ almost every $x \in \mathbb{M}$,
\begin{equation}
\label{eq:noyau}
 P_t f (x) = \int_{\mathbb{M}} f(y) \, p_t(x,dy). 
\end{equation}

Let us consider the following equation
\begin{equation}
 \label{eq:RDP}
\tag*{(\text{$\mathbb{RDP}$})}\left\{
\begin{array}{ll}
 \frac{\partial}{\partial t} \vec{u} (t) & = 
\mathfrak{C} L \vec{u}(t) + G(\vec{u}(t)) \vec{\lambda} , \quad t>0\\
 \vec{u}(0) &= \vec{f} 
\end{array}
\right.
\end{equation}
where, in the two-by-two case,
\begin{itemize} 
 \item the unknown $\vec{u}(t,x)= \left( u_1(t,x), u_2(t,x), u_3(t,x), u_4(t,x) \right)$ 
is a function from $[0,\infty) \times \mathbb{M}$ to $\mathbb{R}^4$; and 
$L \vec{u}= \left( L u_1, L u_2, L u_3, L u_4 \right)$ 
is defined componentwise.
\item $\vec{\lambda} = (\lambda_1, \lambda_2,\lambda_3,\lambda_4) = \lambda (-1, -1,1,1) \in  \mathbb{R}^4$, with $\lambda \in \dR_+$;
\item the nonlinearity $G$ is quadratic: $G(\vec{u}) = u_1 u_2 - u_3 u_4$.
\item $\mathfrak{C}$ is a diagonal matrix of the following form 
$$
\mathfrak{C} = \left( \begin{array}{cccc}
                        C_1 & 0 &0 & 0 \\
                         0 & C_2 & 0 & 0 \\
                         0 & 0 & C_3& 0 \\
                         0 & 0 & 0 & C_4
                       \end{array}
\right),
$$
 where we assume that $C_1=C_3$ and $C_2=C_4$. 

\bigskip 
(This condition is weakened in section \ref{section:extension}). 
\item the initial datum is $\vec{f} = (f_1,f_2,f_3,f_4)$.
\end{itemize}

%

\subsection*{Dirichlet form and logarithmic Sobolev inequality}  
Let $(\mathcal{E}, \mathcal{D})$ be the Dirichlet form associated to $(L, \mu)$ 
(see \cite{davies90}, \cite{fukushima-al}, \cite{ma-rockner}, \cite{bouleau-hirsch}; 
or \cite{FRZ} for a minimal introduction). For any $u \in \mathcal{D}(L)$ 
(the domain of $L$) and
$v \in \mathcal{D}$ (the domain of the Dirichlet form), one has
$$
\mathcal{E}(u,v) = - \mu (v\,Lu  ).
$$  
We will denote $\mathcal E(u)\equiv\mathcal E(u,u)$, for any $u\in\mathcal{D}$. 
Recall that $\mathcal{D}$ is a real Hilbert space with associated norm 
$$
\| u \|_\mathcal{D} = (\mu (u^2) + \mathcal{E}(u))^{1/2}.
$$

We will assume that the Dirichlet structure $(\mathcal E,\mu)$ 
satisfies logarithmic Sobolev inequality 
with constant $C_{LS}\in(0,\infty)$, that is 
\begin{equation}
\label{eq-logsob}
\ent_{\mu}(u^2) \equiv\mu\Big(u^2\log\frac{u^2}{\mu(u^2)}\Big)\leq C_{LS} \mathcal E(u),
\end{equation}
for any $u\in\mathcal D$.


\subsection*{Classical function spaces}

Let $I=[0,T]$.
For any Banach space $(X,\| \cdot \|_X)$, we shall denote
by $C(I,X)$ the Banach space of continuous functions from $I$
to $X$ equipped with the supremum norm 
$$
\sup_{t\in I} \| u(t) \|_X.
$$ 
Let also $\mathbb{L}^2(I,X)$ be the space of 
(a.e. classes of) Bochner measurable functions
from $I$ to $X$
such that
$
\int_0^T \| u(t) \|_X^2 ds<\infty.
$ 
As for vector valued functions, 
let $\mathbb{L}^2(I,X^4)$ 
be the space of Bochner measurable 
functions $t\in I \mapsto (u_1(t),u_2(t),u_3(t), u_4(t)) \in X^4$ such that 
$$
\int_0^T \sum_{i=1}^4 \|u_i(s)\|_X^4 ds<\infty.
$$ 
All these are Banach spaces.

\medskip 
We'll furthermore consider the space $\mathbb{L}^\infty(I,X)$
of Bochner measurable $X$-valued functions on $I$ such that
$$
\ESSUP{0\leq t \leq T} \| u(t) \|_{X} < + \infty.
$$ 

The reader may refer to \cite{schwabik-ye} for Bochner measurability, Bochner integration and other
Banach space integration topics.

\subsection*{Bochner measurability in an Orlicz space}



%

Let $\Phi : \dR \rightarrow \dR_+$ given by $\Phi(x) = \exp(|x|) -1$ 
and $\Phi_\alpha(x) = \Phi(|x|^\alpha)$, $\alpha \geq 1$. 
These  are Young functions and the Orlicz space associated to $\Phi_\alpha$  is 
denoted by $\mathbb{L}^{\Phi_\alpha}(\mu)$. This is the 
space of measurable functions $f$ such that 
\begin{equation}
\label{eq:Phi_alpha_integrable}
\mu( \Phi_\alpha(\gamma f)) < \infty
\end{equation}
for {\it some} $\gamma >0$ (or functions {\it whose $\alpha$ power 
is exponentially integrable}).

An important closed subspace $E^{\Phi_\alpha}(\mu)$ of 
$\mathbb{L}^{\Phi_\alpha}(\mu)$ consists of those
functions such that \eqref{eq:Phi_alpha_integrable} holds for 
{\it any} $\gamma >0$.
This is the closure of the space of simple functions (finitely valued measurable functions)
in $\mathbb{L}^{\Phi_\alpha}(\mu)$.

\medskip
A stricking property of Markov semigroups is that $C^0$ property
in $\mathbb{L}^2(\mu)$
implies $C^0$ property
in any $\mathbb{L}^p(\mu)$; $1 \leq p < +\infty$ (see \cite{davies90}). 
We will need the following weakened result in the context of Orlicz spaces.

\begin{proposition}
\label{proposition:Bochner_measurability_semigroup}
Let $f \in E^{\Phi_\alpha}(\mu)$, $\alpha \geq 1$. Then the linear semigroup is 
Bochner measurable in time in $E^{\Phi_\alpha}(\mu)$. More precisely, 
the mapping 
$t \in [0, \infty) \mapsto P_t f \in E^{\Phi_\alpha}(\mu)$ 
belongs to $\mathbb{L}^\infty([0,\infty), E^{\Phi_\alpha}(\mu))$
and 
$$
\ESSUP{0\leq t < \infty} \| P_t f \|_{E^{\Phi_\alpha}(\mu)} 
\leq \| f \|_{E^{\Phi_\alpha}(\mu)}.
$$
\end{proposition}

The proof is given in appendix \ref{appendix:Bochner-measurability}.


\subsection*{First regularity result and weak solutions}

The following lemma  exhibits the main role
the entropic inequality (see appendix~\ref{appendix_entropic}) 
and the logarithmic Sobolev inequality play to deal with 
the nonlinearity we consider. 
In short, the multiplication 
operator by a function in $\mathbb{L}^{\Phi_2}(\mu)$ 
is a bounded operator, mapping the domain of the Dirichlet form 
$\mathcal{D}$ to $\mathbb{L}^2(\mu)$.

\begin{elem}[Regularity property]
 \label{lemma:regularity_nonlinearity}
Assume the Dirichlet structure $(\mu, \mathcal{E})$ satisfies logarithmic Sobolev
inequality with constant $C_{LS} \in (0, \infty)$. 
Let $\Phi(x) = \exp(|x|)-1$
and ${\Phi_2}(x)= \Phi(x^2)$.
Let $u \in \mathbb{L}^2(I, \mathcal{D}) $
and $v \in 
\mathbb{L}^{\infty}(I, \mathbb{L}^{\Phi_2}(\mu))$.
Then 
$
u v \in \mathbb{L}^2(I,\mathbb{L}^2(\mu)) 
$
and the bilinear mapping 
\begin{equation}
\label{eq:bilinear_continuity_time_spaces}
 (u,v) \in \mathbb{L}^2(I, \mathcal{D}) \times 
\mathbb{L}^{\infty}(I, \mathbb{L}^{\Phi_2}(\mu)) \mapsto u \, v \in \mathbb{L}^2(I,\mathbb{L}^2(\mu))
\end{equation}
is continuous. Consequently,
\begin{equation}
\label{eq:trilinear_continuity_time_spaces}
 (\phi,u,v) \in \mathbb{L}^2(I,\mathbb{L}^2(\mu)) \times \mathbb{L}^2(I, \mathcal{D}) \times 
\mathbb{L}^{\infty}(I, \mathbb{L}^{\Phi_2}(\mu)) \mapsto \phi \, u \, v \in \mathbb{L}^1(I,\mathbb{L}^1(\mu))
\end{equation}
is trilinear continuous.
\end{elem}

The reader may note that
we will use this lemma 
to define  properly a weak solution of the nonlinear problem below. 

\bigskip
\begin{eproof}
Note that $f \in \mathbb{L}^{\Phi_2}(\mu)$ 
iff $f^2 \in \mathbb{L}^\Phi(\mu)$ 
and that 
\begin{equation}
\label{eq:norms}
 \|f^2 \|_\Phi = \| f \|_{\Phi_2}^2.
\end{equation}

\bigskip
First we show that the bilinear mapping
\begin{equation}
\label{eq:bilinear_continuity}
(\su, \sv) \in \mathcal{D} \times \mathbb{L}^{\Phi_2}(\mu) \mapsto \su \sv \in \mathbb{L}^2(\mu)
\end{equation}
is continuous. 
Fix $0<\gamma < \| \sv^2 \|^{-1}_{\mathbb{L}^{\Phi}(\mu))}$.
Then, $\mu(\exp(\gamma \sv^2))-1 \leq 1$ and so $\mu(e^{\gamma \sv^2}) \leq 2$.
Hence, using the entropic inequality~\eqref{eq:entropic_inequality}, and then the logarithmic Sobolev inequality, 
one gets, 
\begin{multline*}
 \mu(\su^2 \sv^2 ) \leq \frac{1}{\gamma} \mu \left( \su^2 \log ( \frac{\su^2}{\mu(\su^2)}) \right)
+ \frac{\mu(\su^2)}{\gamma} \log \mu( e^{\gamma \sv^2}) \\
\leq 
\frac{1}{\gamma} \left( C_{LS} \mathcal{E}(\su) + \log 2 \,  \mu(\su^2) \right) \leq 
\frac{\max(\log 2, C_{LS})}{\gamma} \, \| \su \|_{\mathcal{D}}^2.
\end{multline*}
Letting $\gamma$ go to $\| \sv^2 \|^{-1}_{\mathbb{L}^\Phi(\mu))}$,
using \eqref{eq:norms} one gets the announced continuity. 

If now $u \in \mathbb{L}^2([0,T], \mathcal{D}) $
and $v \in \mathbb{L}^{\infty}([0,T], \mathbb{L}^{\Phi_2}(\mu))$,
there exist two sequences of simple functions (see \cite{schwabik-ye} if necessary) 
$(u_n)_n \subset \mathcal{S}_{I, \mathcal{D}}$ and 
$(v_n)_n \subset \mathcal{S}_{I, \mathbb{L}^{\Phi_2}}$ converging to $u$ (resp. $v$)
a.e. in $\mathcal{D}$ (resp.  $\mathbb{L}^{\Phi_2}$). The continuity of 
\eqref{eq:bilinear_continuity} shows that $(u_n \, v_n)_n$ is a sequence
of simple functions with values in $\mathbb{L}^2(\mu)$ 
which converges a.e. in $\mathbb{L}^2(\mu)$ to $u \, v$. Bochner measurability
of $u \,v$ from $I$ to $\mathbb{L}^2(\mu)$ follows. 

As for continuity of \eqref{eq:bilinear_continuity_time_spaces}, what precedes shows that,  
 for any $t$ a.e.,
\begin{equation*}
\| u(t) v(t) \|_{\mathbb{L}^2(\mu)}^2\leq 
\max(\log 2, C_{LS}) \, \|v\|^2_{\mathbb{L}^\infty([0,T],\mathbb{L}^{\Phi_2}(\mu))} \, 
\| u(t) \|_{\mathcal{D}}^2.
\end{equation*}
Integrating w.r.t. $t$ on $[0,T]$, one gets the result. 
Finally, 
continuity of the trilinear mapping 
follows by Cauchy-Schwarz inequality in $\mathbb{L}^2$.
\end{eproof}

\bigskip
\noindent \textbf{Weak solutions}. Let $T>0$.  We  say that a function 
\begin{equation}
\label{eq:good_spaces} 
\vec{u} \in \bigl(\mathbb{L}^2([0,T], \mathcal{D}) \cap C([0,T], \mathbb{L}^2(\mu)) 
\cap \mathbb{L}^\infty([0,T], \mathbb{L}^{\Phi_2}(\mu)) \bigr)^4
\end{equation}
is a {\it weak solution} of $\ref{eq:RDP}$ on $[0,T]$ provided, for any $\vec{\phi} \in 
{C}^\infty([0,T], \mathcal{D}^4 )$ and any $t \in [0,T)$,

\begin{multline}
\label{eq:weak_RDP}
\tag*{(\text{weak-$\mathbb{RDP}$})}
- \int_0^t \sum_{i=1}^4 \mu ( u_i(s) \partial_s \phi_i (s) ) \, ds  +  \left[  \sum_{i=1}^4 
\mu ( u_i(t) \phi_i (t)  - u_i(0) \phi_i(0))\right] \\
= - \int_0^t  \sum_{i=1}^4  C_i \mathcal{E}(u_i(s), \phi_i(s) ) \, ds 
+  \int_0^t  \sum_{i=1}^4 \lambda_i \mu (\phi_i(s) G(\vec{u}(s)) ) \, ds. 
\end{multline}

When this is satisfied for any $T >0$, we'll say that $\vec{u}$ is a weak solution
on $[0,\infty)$.

\subsection*{Main result}

\begin{ethm}  
\label{thm:main-result}
Let $(L,\mu)$ be a selfadjoint Markov generator satisfying  logarithmic Sobolev inequality~\eqref{eq-logsob}
with constant $C_{LS} \in (0, \infty)$. 

Let $\Phi_2(x) = \exp(x^2) -1$.

Assume $\vec{f} \geq 0$ is a nonnegative initial datum and 
$\vec{f} \in (E^{\Phi_2}(\mu))^4$.  

Then, for any diffusion coefficients $C_1>0$ and $C_2>0$ and any reaction rate $\lambda >0$,  
there exists a unique nonnegative
weak solution $\vec{u}$ of $\ref{eq:RDP}$ on $[0, \infty)$. 

Moreover, for any 
$\alpha \geq 1$, any $\gamma>0$ and any $i=1, \dots, 4$, if $\mu(e^{\gamma (f_1 + f_3)^\alpha})<\infty$ and 
$\mu(e^{\gamma (f_2 + f_4)^\alpha})<\infty$, then 
$$
\forall t \text{ a.e. in } [0, \infty),  \qquad
\mu(e^{\gamma u^\alpha_i(t)}) \leq \max \Bigl(\mu(e^{\gamma (f_1 + f_3)^\alpha}), 
\mu(e^{\gamma (f_2 + f_4)^\alpha}) \Bigr).
$$  
\end{ethm}

In section \ref{section:extension}, we will state  the extension 
of this theorem to the general problem~\eqref{eq-gene}.

\bigskip
In short, to prove this theorem,
we linearize the system of equations by means of an approximation
sequence $(\vec{u}^{(n)})_n$. We show recursively that 
$\vec{u}^{(n)}(t)$ is nonnegative, 
belongs to $\mathbb{L}^\infty([0,T],\mathbb{L}^{\Phi_2}(\mu))$
so that lemma \ref{lemma:regularity_nonlinearity}
guarantees $\vec{u}^{(n+1)}$
is well defined. 
This propagation is made precise in a lemma
studying the linear {\it cornerstone} problem which underlies 
the recursive approach.

We will first focus our efforts to prove convergence of the approximation sequence
in the space $\bigl(\mathbb{L}^2([0,T], \mathcal{D}) 
\cap C([0,T], \mathbb{L}^2(\mu))\bigr)^4$. 
Afterwards, we detail a way to study 
the cornerstone existence lemma. 

\begin{remark}
We will exhibit in appendix~\ref{appendix:Orlicz} a sufficient
condition to ensure that $f \in E^{\Phi_\alpha}(\mu)$, namely, that there exist $\beta >\alpha$ and $\gamma >0$
such that $\mu ( e^{\gamma |f|^\beta} )< +\infty$. In particular, 
it implies that, provided $\vec{f} \geq 0$ belongs to 
$(E^{\Phi_2}(\mu))^4$,
one may choose $\tilde{\gamma}>0$ large enough such that
  \begin{equation}
 \label{eq:constants_constraint}
 \begin{array}{l} 
\disp\frac{4}{\min(C_1,C_2)} \,  \lambda C_{LS} <  \tilde{\gamma}, 
\text{ and } \\
 \vec{f} \text{ still satisfies } \mu(e^{\tilde{\gamma} f_i}) < \infty,
i =1, \dots, 4
\end{array}
\end{equation}
 which will be useful in the proof of existence and uniqueness. 
\end{remark}

\section{Iterative procedure}
\label{section:iterative_approach}
Let us define the approximation sequence 
$(\vec{u}^{(n)})_{n\in\dN}$ in the following way. (First of all, note
the parenthesis in $\vec{u}^{(n)}$ has nothing to do with differentiation,
and has been introduced to distinguish the index from powers).
\begin{itemize}
\item for all $n\in\dN$, $\vec{u}^{(n)}(t=0)=\vec{f} \in (E^{\Phi_2}(\mu))^4$;
\item for $n=0$, $\partial_t \vec{u}^{(0)}(t)=\mathfrak{C}L\vec{u}^{(0)}(t), \quad t >0$;
\item for any $n\geq1$, and $t>0$,
\begin{equation*}
 \left\{
\begin{array}{l}
\partial_t u_1^{(n)}(t)=
C_1 Lu_1^{(n)}(t)- \lambda 
\Bigl(u_2^{(n-1)}(t) u_1^{(n)}(t) - u_4^{(n-1)}(t) u_3^{(n)}(t)\Bigr),\\
\partial_t u_3^{(n)}(t)=
C_1 Lu_3^{(n)}(t)+ \lambda 
\Bigl(u_2^{(n-1)}(t) u_1^{(n)}(t) - u_4^{(n-1)}(t) u_3^{(n)}(t)\Bigr),
\end{array}
\right. \\
\end{equation*}
\begin{equation}
\label{eq-iteration}
\tag*{(\text{$\mathbb{RDP}_n$})}
\text{ ~ } 
\end{equation}
\begin{equation*}
\left\{
\begin{array}{l}
\partial_t u_2^{(n)}(t)=
C_2 Lu_2^{(n)}(t)- \lambda 
\Bigl(u_1^{(n-1)}(t) u_2^{(n)}(t) - u_3^{(n-1)}(t) u_4^{(n)}(t)\Bigr),\\
\partial_t u_4^{(n)}(t)=
C_2 Lu_4^{(n)}(t)+ \lambda 
\Bigl(u_1^{(n-1)}(t) u_2^{(n)}(t) - u_3^{(n-1)}(t) u_4^{(n)}(t)\Bigr).
\end{array}
\right.
\end{equation*}
\end{itemize}
Knowing 
$\vec{u}^{(n-1)} \in \bigl(\mathbb{L}^2([0,T], \mathcal{D}) \cap C([0,T], \mathbb{L}^2(\mu)) 
\cap \mathbb{L}^\infty([0,T], \mathbb{L}^{\Phi_2}(\mu)) \bigr)^4$, 
(which is the case for any $T>0$
under our hypothesis for $\vec{u}^{(0)}$ by proposition 
\ref{proposition:Bochner_measurability_semigroup}), 
this system may be reduced to the four independant
affine scalar equations, with $t>0$,
\begin{equation}
\label{eq:decoupled_system} 
\left\{
\begin{array}{l}
 \partial_t {u}_1^{(n)}={C_1}L{u}_1^{(n)}
 -\lambda  P_{C_2 t} (f_2 +f_4) \, u_1^{(n)} + \lambda {P_{C_1t}(f_1+f_3)} u_4^{(n-1)},\\
\partial_t {u}_3^{(n)}={C_1}L{u}_3^{(n)}
 -\lambda  P_{C_2 t} (f_2 +f_4) \, u_3^{(n)} + \lambda {P_{C_1t}(f_1+f_3)} u_2^{(n-1)},\\
\partial_t {u}_2^{(n)}={C_2}L{u}_2^{(n)}
 -\lambda  P_{C_1 t} (f_1 +f_3) \, u_2^{(n)} + \lambda {P_{C_2t}(f_2+f_4)} u_3^{(n-1)},\\
\partial_t {u}_4^{(n)}={C_2}L{u}_4^{(n)}
 -\lambda  P_{C_1 t} (f_1 +f_3) \, u_4^{(n)} + \lambda {P_{C_2t}(f_2+f_4)} u_1^{(n-1)}.\\
\end{array}
\right.
\end{equation}
whose existence, uniqueness and positivity on $[0,T]$
follows from 
Lemma~\ref{lem-lineaire} below, 
with $A(t)=\lambda P_{C_2 t} (f_2 +f_4)$ and $B(t)=\lambda{P_{C_1t}(f_1+f_3)}u_4^{(n-1)}$
(or similarly), using proposition \ref{proposition:Bochner_measurability_semigroup}
and lemma \ref{lemma:regularity_nonlinearity}.

\begin{elem}[Cornerstone existence lemma]
\label{lem-lineaire}
Let $L$ be a Markov generator satisfying  logarithmic Sobolev inequality
with constant $C_{LS} \in (0,\infty)$. 
Let $T>0$ and $A=A(t) \in \mathbb{L}^\infty([0,T], \mathbb{L}^{\Phi_2}(\mu))$ and 
$B\in \mathbb{L}^2([0,T], \mathbb{L}^2(\mu))$.
Then the Cauchy problem 
\begin{equation}
\label{eq:cornerstone_linear}
\tag*{(\text{$\mathbb{CS}$})}
 \left\{ \begin{array}{l}
          \partial_t u(t) = L u(t)- A(t) \, u(t) + B(t),\\
           u(0)= f, f \in \mathbb{L}^2(\mu)
         \end{array}
\right.
\end{equation}
has a unique weak solution on $[0,T]$. 
Futhermore, 
provided $f$, $A$ and $B$ are assumed nonnegative, 
then the solution $u$ is nonnegative.
\end{elem}

Although classical, we recall that $u \in \mathbb{L}^2([0,T], \mathcal{D})\cap C([0,T],\mathbb L^2(\mu))$
is a weak solution
of \ref{eq:cornerstone_linear}
provided, for any $\phi \in {C}^\infty([0,T], \mathcal{D})$, 
and any $0 \leq t \leq T$, 
\begin{multline}
\label{eq:weak_CS}
\tag*{(\text{weak-$\mathbb{CS}$})}
- \int_0^t  \mu ( u(s) \partial_s \phi (s) ) \, ds  +    
\mu \left( u(t) \phi(t)  - u(0) \phi(0)\right)
= - \int_0^t    \mathcal{E}(u(s), \phi(s) ) \, ds \\
+  \int_0^t  \mu \left(\phi(s) \big[-A(s){u(s)}+B(s)\big] \right) \, ds.
\end{multline}


Recursive equivalence of both systems \ref{eq-iteration} 
and \eqref{eq:decoupled_system}  may be seen as follows. 
Starting from \ref{eq-iteration}, one easily gets  
$$
\left\{
\begin{array}{r}
\displaystyle \partial_t (u_1^{(n)}+u_3^{(n)})=C_1L( u_1^{(n)}+u_3^{(n)})\\
\displaystyle \partial_t (u_2^{(n)}+u_4^{(n)})=C_2L( u_2^{(n)}+u_4^{(n)})
\end{array}
\right.
$$
and writting $u_3^{(n)}(t)=P_{C_1t}(f_1+f_3)-u_1^{(n)}(t)$ (and similarly for the other 
coordinates) gives the announced decoupled system. Conversely, deducing from 
the decoupled system that $u^{(n)}_1+u^{(n)}_3= 
P_{C_1t}(f_1+f_3)$ (and similarly)
follows by induction and uniqueness in lemma \ref{lem-lineaire}.

\bigskip
To be able to define $\vec{u}^{(n+1)}$, and hence prove that the iterative
sequence is well defined, it remains to check that $u^{(n)}_i \in
\mathbb{L}^\infty([0,T], \mathbb{L}^{\Phi_2}(\mu))$, for all 
$i=1, \dots, 4$. This is based on results stated in appendix 
\ref{appendix:Orlicz_semigroup} and can be shown as follows.
\label{page:LinftyLPhi}

\bigskip
We may focus on $u_1^{(n)}(t)$ by symmetry. By positivity of the $u^{(n)}_i$'s 
and constraint 
$u^{(n)}_1+u^{(n)}_3 = P_{C_1 t} (f_1 + f_3)$, the contraction property 
of the semigroup stated in lemma~\ref{lemma:Jensen_Pt_Orlicz}
implies that, for any $\gamma >0$, for any $t$ a.e.,
\begin{equation}
\label{eq:borne_exponentielle_alpha_n}
\mu \Bigl(e^{\gamma (u_1^{(n)}(t))^2} \Bigr)
\leq \mu(e^{\gamma (f_1+f_3)^2}) < +\infty.
\end{equation}
So that, in particular, for any $t \in [0,T]$,
$u_1^{(n)}(t) \in E^{\Phi_2}(\mu)$. Following lemma 
\ref{lemma:LinftyLPhi}, 
what remains to be checked is Bochner 
measurability of the mapping $t \mapsto u_1^{(n)}(t) \in E^{\Phi_2}(\mu)$.

From the corresponding weak 
formulation \ref{eq:weak_CS} applied to a constant
(in time) test function $\phi(t) \equiv \varphi \in \mathcal{D}$,
\begin{multline*}
\mu( u_1^{(n)}(t) \varphi) 
= \mu( f_1 \varphi) - C_1 
\int_0^t \mathcal{E}(u_1^{(n)}(s), \varphi) ds \\
+ \int_0^t \mu \Bigl( \varphi (-\lambda  P_{C_2 s} (f_2 +f_4) \, u_1^{(n)}(s) + 
\lambda {P_{C_1s}(f_1+f_3)} u_4^{(n-1)}(s) \Bigr) ds.
\end{multline*}
Hence, the function $t \mapsto \mu( u_1^{(n)}(t) \varphi)$ is continuous,
for any fixed $\varphi \in \mathcal{D}$. Now, 
$\mathcal{D}$ is a dense subspace of the dual space 
$(E^{\Phi_2})'= \mathbb{L}^{\Phi_2^\ast}(\mu)$ (see appendix
\ref{appendix:Orlicz_semigroup}), so that weak measurability of  $t \mapsto 
u_1^{(n)}(t) \in E^{\Phi_2}$ follows. 
By Pettis measurability theorem\footnote{\label{pettis} see \cite{yosida}, 
\cite{diestel-uhl} or 
\cite{schwabik-ye} for a proof,
\cite{evans}, appendix E.5, theorem 7, for a statement.}
and separability of
$E^{\Phi_2}(\mu)$,  
$t \in [0,T] \mapsto u_1^{(\infty)}(t) \in E^{\Phi_2}(\mu)$ is 
Bochner measurable.

%

\section{Proof of Theorem~\ref{thm:main-result}}

\label{section:proof_main}

\subsection*{Convergence of the approximation procedure (\text{$\mathbb{RDP}_n$})}
From now on, we'll use the notation 
$$
|\vec{u}|^2 =\sum_{i=1}^4 u_i^2 \quad \text{ and } \quad \mathcal E(\vec{u})=\sum_{i=1}^4 \mathcal E({u}_i).
$$
The main idea is to show that, with 
\begin{equation}
\label{eq-def-sigma}
\Sigma_n(t)=\mu(|\vec{u}^{(n)}-\vec{u}^{(n-1)}|^2)(t)+ 2 \kappa\int_0^t\mathcal E(\vec{u}^{(n)}-\vec{u}^{(n-1)})(s)ds,
\end{equation}
for some $\kappa>0$ (specified later), 
the supremum $\sup_{t\in[0,T]} \Sigma_n(t)$ goes to 0 exponentially fast 
as $n$ goes to $\infty$ provided $T>0$ is small
enough. 
From lemma \ref{lem-lineaire},
$\vec{u}^{(n)}$ is defined recursively as a weak solution of 
the cornerstone linear problem.
To make things simpler at this stage, we here perform formal computations
to get a priori estimates. Getting the estimates rigorously makes use of
Steklov regularisation, 
which we will illustrate in the proof of the next proposition.

\subsubsection*{Estimate of the $\mathbb{L}^2$-norm derivative}
We will focus on the $\mathbb{L}^2$-norm of  $u_1^{(n)}$. 
\begin{multline*}
\frac12\frac{d}{dt}\mu\big[ (u_1^{(n)}-u_1^{(n-1)})^2\big]=C_1\mu\big[(u_1^{(n)}-u_1^{(n-1)})L(u_1^{(n)}-u_1^{(n-1)})\big]\\
\quad-\lambda\mu\big[(u_1^{(n)}-u_1^{(n-1)})(u_1^{(n)}u_2^{(n-1)}-u_3^{(n)}u_4^{(n-1)}-u_1^{(n-1)}u_2^{(n-2)}+u_3^{(n-1)}u_4^{(n-2)})\big],
\end{multline*}
and after natural multilinear handlings,
\begin{multline*}
\frac12\frac{d}{dt}\mu\big[ (u_1^{(n)}-u_1^{(n-1)})^2\big]=-C_1\mathcal E\big[u_1^{(n)}-u_1^{(n-1)}\big]-\lambda\mu\big[(u_1^{(n)}-u_1^{(n-1)})^2u_2^{(n-1)}\big]\\
-\lambda\mu\big[(u_1^{(n)}-u_1^{(n-1)})(u_2^{(n-1)}-u_2^{(n-2)})u_1^{(n-1)}\big]+
\lambda\mu\big[(u_1^{(n)}-u_1^{(n-1)})(u_3^{(n)}-u_3^{(n-1)})u_4^{(n-1)}\big]\\
+\lambda\mu\big[(u_1^{(n)}-u_1^{(n-1)})(u_4^{(n-1)}-u_4^{(n-2)})u_3^{(n-1)}\big].
\end{multline*}
Since $\vec{u}^{(n-1)}$ is nonnegative, using the quadratic inequality  
$ab\leq a^2/2+b^2/2$, 
one gets
\begin{multline*}
\frac12\frac{d}{dt}\mu\big[ (u_1^{(n)}-u_1^{(n-1)})^2) \leq -C_1\mathcal E\big[u_1^{(n)}-u_1^{(n-1)}\big]\\
+\frac{\lambda}{2}\mu\big[(u_1^{(n)}-u_1^{(n-1)})^2u_1^{(n-1)}\big]+\frac{\lambda}{2}\mu\big[(u_2^{(n-1)}-u_2^{(n-2)})^2u_1^{(n-1)}\big]\\
+\frac{\lambda}{2}\mu\big[(u_1^{(n)}-u_1^{(n-1)})^2u_4^{(n-1)}\big]+\frac{\lambda}{2}\mu\big[(u_3^{(n)}-u_3^{(n-1)})^2u_4^{(n-1)}\big]\\
+\frac{\lambda}{2}\mu\big[(u_1^{(n)}-u_1^{(n-1)})^2u_3^{(n-1)}\big]+\frac{\lambda}{2}\mu\big[(u_4^{(n-1)}-u_4^{(n-2)})^2u_3^{(n-1)}\big].
\end{multline*}
All the similar terms are then estimated  
thanks to the relative entropy inequality~\eqref{eq:entropic_inequality}. 
For
instance, 
\begin{multline*}
\mu\big[(u_1^{(n)}-u_1^{(n-1)})^2u_1^{(n-1)}\big]\leq  \frac{1}{\gamma}\ent_\mu\big[(u_1^{(n)}-u_1^{(n-1)})^2\big]\\
+\frac{1}{\gamma}\mu\big[\big(u_1^{(n)}-u_1^{(n-1)}\big)^2\big]\log\mu\big[e^{\gamma u_1^{(n-1)}}\big].
\end{multline*}
The logarithmic Sobolev inequality~\eqref{eq-logsob} and 
bound \eqref{eq:borne_exponentielle_alpha_n} give 
$$
\mu\big[(u_1^{(n)}-u_1^{(n-1)})^2u_1^{(n-1)}\big]\leq  
\frac{C_{LS}}{\gamma}\mathcal E\big[u_1^{(n)}-u_1^{(n-1)}\big]+
\frac{D}{\gamma}\mu\big[(u_1^{(n)}-u_1^{(n-1)})^2\big],
$$
where 
\begin{equation}
\label{eq-def-D}
D=\max\BRA{\log\mu(e^{\gamma (f_1+f_3)}),\log\mu(e^{\gamma (f_2+f_4)})}.
\end{equation}
Using the same arguments for all the terms leads to  
\begin{multline*}
\frac12\frac{d}{dt}\mu\big[ (u_1^{(n)}-u_1^{(n-1)})^2\big]\leq-C_1\mathcal E\big[u_1^{(n)}-u_1^{(n-1)}\big]\\
+\frac{\lambda C_{LS}}{2\gamma}\Big(3\mathcal E\big[u_1^{(n)}-u_1^{(n-1)}\big]+\mathcal E\big[u_2^{(n-1)}-u_2^{(n-2)}\big]\\
+\mathcal E\big[u_3^{(n)}-u_3^{(n-1)}\big]+\mathcal E\big[u_4^{(n-1)}-u_4^{(n-2)}\big]\Big)\\
+D \, \frac{\lambda}{2\gamma}\Big(3\mu\big[(u_1^{(n)}-u_1^{(n-1)})^2\big]+\mu\big[(u_2^{(n-1)}-u_2^{(n-2)})^2\big]\\
+\mu\big[(u_3^{(n)}-u_3^{(n-1)})^2\big]+\mu\big[(u_4^{(n-1)}-u_4^{(n-2)})^2\big]\Big).
\end{multline*}
Completely similar terms are obtained when dealing with the $\mathbb{L}^2$-norms of the other
components. 
After summation in all the components, one gets
 \begin{multline*}
\frac12\frac{d}{dt}\mu\big[ |\vec{u}^{(n)}-\vec{u}^{(n-1)}|^2\big]
\leq -\min (C_1,C_2) \mathcal E\big[\vec{u}^{(n)}-\vec{u}^{(n-1)}\big]\\
 + \frac{\lambda C_{LS}}{2\gamma} \, \Bigl( 4 \mathcal E\big[\vec{u}^{(n)}-\vec{u}^{(n-1)}\big] 
+ 2 \mathcal E\big[\vec{u}^{(n-1)}-\vec{u}^{(n-2)}\big] \Bigr) \\
\qquad 
+ \frac{D \lambda}{2\gamma} 
\Bigl( 4 \mu\big[ |\vec{u}^{(n)}-\vec{u}^{(n-1)}|^2\big] + 2 \mu\big[ |\vec{u}^{(n-1)}-\vec{u}^{(n-2)}|^2\big] 
\Bigr).
\end{multline*}
Let $\kappa \equiv \min (C_1,C_2) -  \frac{2 \lambda C_{LS}}{\gamma}$
which is positive thanks to the assumed constraint~\eqref{eq:constants_constraint}.

\bigskip
Use the absolute continuity and the positivity of $\int_0^t \mathcal E(\vec{u}^{(n)}-\vec{u}^{(n-1)})(s)ds$, to get
\begin{multline*}
\frac12\frac{d}{dt}\PAR{\mu\big[ |\vec{u}^{(n)}-\vec{u}^{(n-1)}|^2\big]+2\kappa\int_0^t \mathcal E\big[\vec{u}^{(n)}-\vec{u}^{(n-1)}\big](s)ds} \leq\\  
 D \, \frac{2 \lambda}{\gamma}\PAR{ \mu\big[ |\vec{u}^{(n)}-\vec{u}^{(n-1)}|^2\big]+2\kappa\int_0^t \mathcal E\big[\vec{u}^{(n)}-\vec{u}^{(n-1)}\big](s)ds}\\
+ D \, \frac{ \lambda}{\gamma} \PAR{\mu\big[ |\vec{u}^{(n-1)}-\vec{u}^{(n-2)}|^2\big]+\frac{ C_{LS}}{D} \, \mathcal E\big[\vec{u}^{(n-1)}-\vec{u}^{(n-2)}\big]}.
\end{multline*}
Reminding the definition~\eqref{eq-def-sigma} of $\Sigma_n$ and that 
$\vec{u}^{(n)}(0) = \vec{u}^{(n-1)}(0)$, after integration over $[0,t]$,
$t\in[0,T]$, we obtain the following main estimate
\begin{multline}
\label{eq-estimation-l2}
\Sigma_{n}(t)\leq D \, \frac{4 \lambda}{\gamma}\int_0^t\Sigma_{n}(s)ds\\
+D \, \frac{2\lambda}{\gamma} \PAR{\int_0^t\mu\big[ |\vec{u}^{(n-1)}-\vec{u}^{(n-2)}|^2\big](s)ds
+\frac{ C_{LS}}{D} \, \int_0^t\mathcal E\big[\vec{u}^{(n-1)}-\vec{u}^{(n-2)}\big](s)ds}.
\end{multline}
\subsubsection*{Gronwall argument and convergence}
Gronwall type arguments applied to the estimate~\eqref{eq-estimation-l2} give for any $t\in[0,T]$,
\begin{multline*}
\Sigma_n(t)\leq D \, \frac{2 \lambda}{\gamma} e^{D \, \frac{4 \lambda}{\gamma}t}\Big(\int_0^t\mu( |\vec{u}^{(n-1)}-\vec{u}^{(n-2)}|^2)(s)ds\\
+\frac{ C_{LS}}{D} \, \int_0^t\mathcal E(\vec{u}^{(n-1)}-\vec{u}^{(n-2)})(s)ds\Big),
\end{multline*}
It follows that
$$
\sup_{t\in[0,T]}\Sigma_n(t)\leq \eta(T)\sup_{t\in[0,T]}\Sigma_{n-1}(t),
$$
where $\eta(T)= \frac{2\lambda}{\gamma}e^{D \, \frac{4 \lambda}{\gamma}T}(D  T+{ C_{LS}} )$.

Condition~\eqref{eq:constants_constraint} implies that there exists $T(D)>0$ such that 
$\eta(T)<1$ since $\lim_{T\rightarrow0}\eta(T)=\frac{4 \lambda}{\gamma} { C_{LS}}$. Therefore, for this choice of $T>0$, $(\vec{u}^{(n)})_{n\in\dN}$  satisfies 
\begin{multline*}
 \max\BRA{\int_0^T\mathcal E(\vec{u}^{(n)}-\vec{u}^{(n-1)})(s)ds,
\sup_{t\in[0,T]}\mu( |\vec{u}^{(n)}-\vec{u}^{(n-1)}|^2)(s)}\\
 \leq\eta(T)^{n-1}\sup_{t\in[0,T]}\Sigma_1(t),
\end{multline*}
Performing a similar estimate for 
$\frac{1}{2} \frac{d}{dt} \mu( (u_i^{(n)}(t))^2 )$, one gets the uniform bound 
\begin{equation}
\label{eq:uniform_bound_n}
\forall n, \forall t \in [0,T],  \Sigma_n(t) \leq e^{\frac{4 \lambda D}{\gamma} t} \, \mu(|\vec{f}|^2).
\end{equation}
It follows that 
$$
\| \vec{u}^{(n)}-\vec{u}^{(n-1)} \|_{\mathbb{L}^2([0,T],\mathcal D^4)\cap\mathcal C([0,T],\mathbb{L}^2(\mu)^4)}
\leq 2 e^{\frac{4 \lambda D}{\gamma} T} \, \mu(|\vec{f}|^2) \, \eta(T)^{n-1}.
$$
Hence, $(\vec{u}^{(n)})_{n\in\dN}$ is a Cauchy sequence: it converges 
to  some function $\vec{u}^{(\infty)}
\in\mathbb{L}^2([0,T],\mathcal D^4)\cap\mathcal C([0,T],\mathbb{L}^2(\mu)^4)$. 
 
\subsubsection*{Global existence of the weak solution}
\label{weak-solution}
Let $T>0$ fixed as in the previous computation. We will 
first prove that the limit $\vec{u}^{(\infty)}$ 
is a weak solution of~\ref{eq:RDP} in $[0,T]$. 
Let $\phi\in\mathcal C^{\infty}([0,T],\mathcal D)$ 
and use the weak formulation of ~\ref{eq-iteration} for 
$\vec{\phi} \equiv (\phi,0,0,0)$. 
For any $t\in[0,T]$,  
\begin{multline*}
-\int_0^t \mu(  u_1^{(n)}(s) \partial_s\phi ) ds+ \mu(u_1^{(n)}(t)\phi(t) )-\mu(f_1
 \phi(0))=\\
-C_1\int_0^t\mathcal E(\phi, u_1^{(n)})(s)ds
-\lambda\int_0^t\mu( \phi u_1^{(n)} u_2^{(n-1)})(s)ds+ \lambda\int_0^t\mu(\phi u_3^{(n)} u_4^{(n-1)})(s)ds.
\end{multline*}
We now show we can pass to the limit $n \rightarrow \infty$ in all the terms.
(Dealing with other coordinates $u_i^{(n)}$ is similar by symmetry).
Thanks to  the continuity of the  scalar product in $\mathbb L^2([0,T],\mathcal D)$, we have
$$
\lim_{n\rightarrow\infty}\int_0^t \mu(  u_1^{(n)}\partial_s\phi )(s)ds=\int_0^t \mu(  u_1^{(\infty)}\partial_s\phi )(s)ds
$$
and 
$$
\lim_{n\rightarrow\infty}\int_0^t\mathcal E(\phi, u_1^{(n)})(s)ds=\int_0^t\mathcal E(\phi, u_1^{(\infty)})(s)ds.
$$
Moreover, as the convergence also holds in $\mathcal C([0,T],\mathbb L^2(\mu))$, then  
$
\lim_{n\rightarrow\infty}\mu(  u_1^{(n)}\phi )(t)=\mu(  u_1^{(\infty)}\phi )(t)
$
and
$\mu( f_1 \phi(0) ) = \lim_{n\rightarrow\infty}\mu(  u_1^{(n)} \phi )(0) =
\mu(  u_1^{(\infty)}\phi )(0).
$

\bigskip
Dealing with the convergence of the term $\int_0^t\mu( \phi u_1^{(n)}u_2^{(n-1)})(s)ds$ 
(and similarly of $\int_0^t\mu( \phi u_3^{(n)}u_4^{(n-1)})(s)ds$) 
is more  intricate. The difficulty is to show that $u^{(\infty)}_1$ belongs to 
$\mathbb{L}^\infty([0,T],E^{\Phi_2})$ which will follow indirectly. 
The details are as follows.

By lemma \ref{lemma:regularity_nonlinearity}, $\tau_n \equiv \tau^{(12)}_n \equiv  
\mu( \phi u_1^{(n)}u_2^{(n-1)}) \in \mathbb{L}^1([0,T])$. Let us show that
this sequence is Cauchy, and so converges to, say, $\tau^{(12)}$ in $\mathbb{L}^1([0,T])$.
Indeed,
\begin{multline*}
\|\tau_n - \tau_m \|_1 \leq  
\int_0^t\mu( |\phi(s)| \cdot |u_1^{(n)}u_2^{(n-1)}(s) -u_1^{(m)}u_2^{(m-1)}(s)|) ds
\leq\\
\int_0^t\mu( |\phi| \cdot | u_1^{(n)} -u_1^{(m)}| \cdot |u_2^{(n-1)}|)(s)ds
+\int_0^t\mu( |\phi| \cdot | u_2^{(n-1)} -u_2^{(m-1)}| \cdot |u_1^{(m)}|)(s)ds.
\end{multline*}
But by \eqref{eq:borne_exponentielle_alpha_n}, and 
again entropic and log-Sobolev inequalities,
\begin{multline}
 \label{eq:nonlinear_term_convergence}
\int_0^t\mu( |\phi| \cdot | u_1^{(n)} -u_1^{(m)}| \cdot |u_2^{(n-1)}|)(s)ds \\
\leq \frac{\max(C_{LS}, \log M_\gamma)}{\gamma} 
\, \| \phi \|_{\mathbb{L}^2(I, \mathbb{L}^2(\mu))} \, 
\| u_1^{(n)} -u_1^{(m)} \|_{\mathbb{L}^2(I, \mathcal{D})}
\end{multline}
with $M_\gamma \equiv 
\max \Bigl( \mu(e^{\gamma (f_1+f_3)^2}),\mu(e^{\gamma (f_2+f_4)^2})\Bigr)$.
This goes to $0$ as $n,m \rightarrow +\infty$.

Now, for any $t \in [0,T]$, $u_1^{(n)}(t) \rightarrow  u_1^{(\infty)}(t)$ in 
$\mathbb{L}^2(\mu)$, so that along a subsequence it converges $\mu$ a.s..
Hence first $u_1^{(\infty)}(t)$ is nonnegative ($\mu$ a.s.) and secondly 
by Fatou lemma 
$$
\mu(e^{\gamma (u_1^{(\infty)}(t))^2}) \leq \liminf_n \mu(e^{\gamma (u_1^{(n)}(t))^2})
\leq M_\gamma < \infty
$$
for any $t$ a.e. in $[0,T]$. And this for any $\gamma>0$.    
Consequently, for any $t$ a.e., $u^{(\infty)}(t) \in E^{\Phi_2}$.

\bigskip
From lemma \ref{lemma:LinftyLPhi}, what remains to do is to prove 
$E^{\Phi_2}$ Bochner measurability.
Let us summarize what we obtained. One has after taking limit $n \rightarrow +\infty$,
\begin{multline*}
-\int_0^t \mu(  u_1^{(\infty)}(s) \partial_s\phi ) ds+ \mu(u_1^{(\infty)}(t)\phi(t) )-\mu(f_1
 \phi(0))=\\
-C_1\int_0^t\mathcal E(\phi, u_1^{(\infty)})(s)ds
-\lambda\int_0^t\tau^{(12)}(s)ds+ \lambda\int_0^t\tau^{(34)}(s)ds.
\end{multline*}
In particular, choosing $\phi(t) = \varphi \in \mathcal{D}$, 
the mapping $t \in [0,T] \mapsto \mu( \varphi  u_1^{(\infty)}(t)) \in \dR$
is continuous. Then, arguments detailed on page \pageref{page:LinftyLPhi} ensure 
that $u_1^{(\infty)} \in \mathbb{L}^\infty(I,L^{\Phi_2})$. 

Furthermore,
\begin{multline*}
\Big|\int_0^t\mu( \phi u_1^{(n)}u_2^{(n-1)})(s)ds-\int_0^t\mu( \phi u_1^{(\infty)}u_2^{(\infty)})(s)ds\Big|\leq\\
\int_0^t\mu( |\phi| \cdot | u_1^{(n)} -u_1^{(\infty)}|
\cdot |u_2^{(n-1)}|)(s)ds
+\int_0^t\mu( |\phi| \cdot | u_2^{(n-1)} -u_2^{(\infty)}|
\cdot |u_1^{(\infty)}|)(s)ds.
\end{multline*}
Letting separately $n$ (resp. $m$) to $+\infty$ in \eqref{eq:nonlinear_term_convergence}
shows that
$$
\int_0^t\mu( \phi u_1^{(n)}u_2^{(n-1)})(s)ds \rightarrow \int_0^t\mu( \phi u_1^{(\infty)}u_2^{(\infty)})(s)ds.
$$   
All this implies that 
$$
\vec{u}^{(\infty)}=(u_1^{(\infty)},u_2^{(\infty)},u_3^{(\infty)},u_4^{(\infty)})
\in \Bigl(\mathbb{L}^2(I,\mathcal D)\cap  C(I,\mathbb{L}^2(\mu)) \cap 
\mathbb{L}^\infty(I,E^{\Phi_2}) \Bigr)^4
$$ 
is a nonnegative weak solution of~\ref{eq:RDP}.

From the local existence in $[0,T]$ 
to a global existence in $[0,\infty)$ 
it is enough to prove that we can repeat the method on the interval 
$[T,2T]$. This follows from the estimate 
$$
\mu(e^{\gamma (u_1^{(\infty)}+u_3^{(\infty)})^2(T)})\leq\mu(e^{\gamma (f_1+f_3)^2}).
$$
See lemma \ref{lemma:Jensen_Pt_Orlicz}.
%

\begin{eprop}[Uniqueness]
\label{prop-unique}
Let $\vec{f} \geq 0$ such that, for some $\gamma >0$, 
$$
M \equiv \max\left\{ \mu( e^{\gamma (f_1 + f_3)}), \mu(e^{\gamma (f_2 + f_4)})\right\} < \infty. 
$$
Assume the diffusion coefficients $C_1$ and $C_2$, the logarithmic Sobolev constant $C_{LS}$ of $L$,
the reaction rate $\lambda$ and the exponential integrability parameter $\gamma$ 
are linked by the constraint
$$
4 \frac{ \lambda C_{LS}}{\min(C_1,C_2)} \leq \gamma.
$$
Then a weak solution of the Reaction-Diffusion problem \ref{eq:RDP} with initial datum $\vec{f}$
is unique.
\end{eprop}

\bigskip We recall basics on Steklov calculus (see \cite{ladyzenskaja} for instance), 
i.e. appropriate time regularization
to deal with weak solutions. For any Banach space
$X$, and any $v \in \mathbb{L}^2([0,T], X)$, the Steklov average\label{steklov}, defined by
\begin{equation*}
a_{h}(v) (t) =  \left\{ \begin{array}{ll} \frac{1}{h} \int_t^{t+h}
v(\tau) \, d \tau &, 0 \leq t \leq T-h,\\
 0 &, T-h < t \leq T
 \end{array}
 \right.
\end{equation*}
converges to $v$ in $\mathbb{L}^2([0,T], X)$ when $h$ goes to $0$.
Moreover, provided $v \in C([0,T],X)$, $a_{h}(v) \in
C^1([0,T-h],X)$, $ \frac{d}{dt} a_h(v) (t) = \frac{1}{h} (v(t+h) - v(t))$
in $X$, and $a_{h}(v)(t)$ converges to $v(t)$ in $X$,  for every $t\in[0,T]$. The space $X$
will be here $\mathbb{L}^2(\mu)$ or $\mathcal{D}$
depending on the context.

\bigskip

 {\noindent {{\it \textbf{Proof of Proposition~\ref{prop-unique}}}}\\\proofbegin~}
Let $\vec{u}$ and $\vec{v}$ be two weak solutions of \ref{eq:RDP} with the same initial datum 
$\vec{f} \geq 0$.
Let $M \in (0,\infty)$ such that, $\forall i =1, \dots, 4$, 
$\mu(e^{\gamma |u_i(t)|}) \leq M$, $t$ a.e., (and similarly for $\vec{v}$). 
Let $\vec{w} \equiv \vec{u} - \vec{v}$ 
and $a_h(w^i) (t)$ the Steklov average of the $i$ component of $w$ 
as defined before. Integrating 
$\frac{1}{2} \frac{d}{ds} \mu \Bigl( (a_h(w^i) (s))^2 \Bigr) 
= \mu (a_h(w^i) (s) \partial_s a_h(w^i) )$ one gets
\begin{equation}
\label{eq:fundamental_thm_Steklov_uniqueness}
\mu \Bigl( (a_h(w^i) (t))^2 \Bigr) = \mu \Bigl( (a_h(w^i) (0))^2 \Bigr)
+ 2 \int_0^t ds \mu \bigl(a_h(w^i) (s) \frac{1}{h} (w^i(s+h) - w^i(s))\bigr).
\end{equation}
We then use the definition of a weak solution with the constant test function
$a_h(w^i) (s) \in  \mathcal{D}$  on the interval $[s,s+h]$
to get 
\begin{multline*}
 \mu( a_h(w^i) (s) \frac{1}{h} (w^i(s+h) - w^i(s))) = -C_i \frac{1}{h} 
\int_s^{s+h} \mathcal{E}(a_h(w^i) (s), w^i(\tau)) d\tau \\
+ \lambda_i \frac{1}{h} 
\int_s^{s+h} d\tau \mu \Bigl(a_h(w^i) (s) \bigl\{ (u_1 u_2 - u_3 u_4) (\tau) - (v_1 v_2 - v_3 v_4) (\tau)
\bigr\} \Bigr)
\end{multline*}
Now, first,
$$
\frac{1}{h} 
\int_s^{s+h} \mathcal{E}(a_h(w^i) (s), w^i(\tau)) d\tau =  \mathcal{E}(a_h(w^i) (s), a_h(w^i) (s)).
$$
And the other term is bounded from above by
\begin{multline*}
 \frac{\lambda}{h^2} \int_{[s, s+h]^2}  d\tau d\tau' \mu \Bigl(|u_i-v_i|(\tau') \,  
\bigl\{ |u_1 - v_1| (\tau) |u_2| (\tau) + |v_1|(\tau) |u_2 -v_2|(\tau) +  \\
|u_3 - v_3| (\tau) |u_4| (\tau) + |v_3|(\tau) |u_4 -v_4|(\tau)
\bigr\} \Bigr).
\end{multline*}
We can deal with the four similar terms by the same way: let us focus
on the first one. One first uses 
$$
\mu \Bigl(|u_i-v_i|(\tau') 
|u_1 - v_1| (\tau) |u_2| (\tau) \Bigr) \leq \frac{1}{2} \mu \Bigl((u_i-v_i)^2(\tau') |u_2| (\tau)
+ (u_1 - v_1)^2 (\tau) |u_2| (\tau) \Bigr).
$$ 
Once gain, entropic inequality followed
by logarithmic Sobolev inequality give
\begin{multline*}
 \frac{\lambda}{h^2} \int_{[s, s+h]^2}  d\tau d\tau' \mu \Bigl(|u_i-v_i|(\tau') 
|u_1 - v_1| (\tau) |u_2| (\tau) \Bigr) \leq \\
\leq \frac{\lambda}{2 \gamma} \frac{1}{h}  
\int_s^{s+h} \Bigl( C_{LS} \mathcal{E}(u_i-v_i)(\tau')  
+ {\log M} \mu( (u_i-v_i)^2(\tau')) \Bigr) d\tau' \\
+ \frac{\lambda}{2 \gamma} \frac{1}{h}  
\int_s^{s+h} \Bigl( C_{LS} \mathcal{E}(u_1-v_1)(\tau)  
+ {\log M} \mu( (u_1-v_1)^2(\tau)) \Bigr) d\tau.
\end{multline*}
Note that, up to a constant, the first term of the RHS is the Steklov average
of the $\mathbb{L}^1([0,T])$ function $ C_{LS} \mathcal{E}(u_i-v_i)(\cdot)  
+ {\log M} \mu( (u_i-v_i)^2(\cdot))$, so that, as $h \rightarrow 0$, 
it converges in $\mathbb{L}^1([0,T])$
to that function. Going back to \eqref{eq:fundamental_thm_Steklov_uniqueness}
and performing all the explained bounds before passing to the limit $h \rightarrow 0$, one
gets the estimate (note that $w_i(0)=0$)
\begin{multline*}
 \mu (w_i^2(t)) \\
\leq 2 \, \int_0^t \! \! \! ds \Bigl( -C_i \mathcal{E}(w_i)(s) + 
\frac{\lambda C_{LS}}{2 \gamma} \Bigl[ 4 \mathcal{E}(w_i)(s) + \mathcal{E}(\vec{w})(s)\Bigr] 
+  \frac{\lambda \log M}{2 \gamma} \Bigl[ 4 \mu( w_i^2(s)) + \mu( \vec{w}^2)(s) \Bigr] \Bigr)
\end{multline*}
Summing over all $i$'s, one gets
\begin{multline*}
 \mu( \vec{w}^2)(t) \leq 2 \, (-\min(C_1,C_2) + 4 \frac{\lambda C_{LS}}{\gamma}) \int_0^t  ds \mathcal{E}(\vec{w})(s)
+ 8 \frac{\lambda \log M}{\gamma}  \int_0^t  ds \mu( \vec{w}^2)(s) \\
\leq 8 \frac{\lambda \log M}{\gamma}  \int_0^t  ds \mu( \vec{w}^2)(s)
\end{multline*}
provided the announced constraint $ 4 \frac{\lambda C_{LS}}{\gamma} \leq \min(C_1,C_2)$
is satisfied. Uniqueness follows by Gronwall arguments.
               {\proofend\\}

\section{Proof of  Lemma~\ref{lem-lineaire}}

\label{section:proof_cornerstone}

Our approach to study the cornerstone linear problem introduced in lemma \ref{lem-lineaire}
will be as follows. We first complete regularity lemma \ref{lemma:regularity_nonlinearity}
by another preliminary lemma (relative to differentiability)
which allow us to perform a recursive approximation of the solution of a mollified 
problem (with a small action of the semigroup on the extra affine term). 
On the way, we 
show a priori estimates which will be useful later to remove the mollification and get
a solution of our initial problem.   
Uniqueness and preservation of positivity
are tackled in specific sections.

Such an approach was already proposed in \cite{FRZ},
and computations look quite similar. The main difference consists in the fact
that, as $A(t) \in \mathbb{L}^{\Phi_2}(\mu)$, then one has 
$\mu(e^{\gamma |A(t)|})<\infty$ for any $\gamma$ (see appendix \ref{appendix:Orlicz}), 
so that, using of the entropic inequality, contribution of the affine extra term may 
be made small enough to be dominated by the log-Sobolev constant without further constraint.

\subsection{Preliminaries}

We recall that $\mathbb{L}^2(\mu)$ may be continuously embedded
in the dual space $\mathcal{D}'$ of the domain $\mathcal{D}$.
From lemma \ref{lemma:regularity_nonlinearity}, it follows that the multiplication
operator by a function $v \in \mathbb{L}^{\infty}([0,T], \mathbb{L}^{\Phi_2}(\mu))$
is a particular case of a Lipschitz continuous operator from
$\mathbb{L}^2([0,T], \mathcal{D})$ to $\mathbb{L}^2([0,T], \mathcal{D}')$.
The following lemma may be stated in this more general context (an example
of which was studied in \cite{FRZ}).

%
%

\begin{elem}[Absolute continuity, differentiability a.e. and weak solutions]
\label{lemma_plus-hypothetic}                      
Let $z \in \mathbb{L}^2([0,T], \mathcal{D}')$, $f \in \mathbb{L}^2(\mu)$ and $\varepsilon >0$.
Define $u(t) = P_t f + \int_0^t P_{t-s+ \varepsilon} z(s) ds$.
Then $u$ belongs to $\mathbb{L}^2([0,T], \mathcal{D}) \cap C([0,T], \mathbb{L}^2(\mu))$ and is (strongly) 
absolutely continuous from  $[a,T]$ to $\mathbb L^2(\mu)$, for any $0<a<T$. And consequently,  the continuous 
function
$
t \in [0,T] \mapsto \mu (u^2(t)) \in \dR
$ 
is absolutely continuous on $[a,T]$.  

Moreover, for all $t$ a.e. in $[0,T]$,
$u(t)$ is differentiable w.r.t. $t$ in $\mathbb{L}^2(\mu)$, belongs to the domain of $L$,
and satisfies 
\begin{equation}
\label{eq:equation_generale_decorelee}
 \left\{ \begin{array}{l}
   \frac{\partial}{\partial t} u(t) = L u(t) + P_\varepsilon (z(t)), \, t \text{ a.e.}\\
   u(0) = f.      
         \end{array}
\right.
\end{equation}
As a consequence, $u \in \mathbb{L}^2([0,T], \mathcal{D}) \cap C([0,T], \mathbb{L}^2(\mu))$ is a weak solution of \eqref{eq:equation_generale_decorelee} i.e., for any $\phi \in\mathcal C^\infty([0,T],\mathcal{D})$,
\begin{multline}
\label{eq:weak_free_general_equation}
- \int_0^t  \mu ( u(s) \partial_s \phi (s) ) \, ds  +    
\mu \left( u(t) \phi(t)  - u(0) \phi(0)\right)
= - \int_0^t    \mathcal{E}(u(s), \phi(s) ) \, ds \\
+  \int_0^t  \mu \left(\phi(s) P_{\varepsilon} z(s) \right) \, ds. 
\end{multline}
\end{elem}

\begin{eproof}
Let us first note that  the Markov semigroup itself satisfies all the announced 
assertions. We only focus on absolute continuity.

Let $([a_i,b_i])_{i=1, \dots, N}$ be a finite collection of (non empty) non overlapping 
subintervals of $[a,T]$. Then, 
$
P_{b_i} f - P_{a_i} f = \int_0^{b_i-a_i} P_\tau L P_{a_i} f d\tau
$
so that 
$$
\| P_{b_i} f - P_{a_i} f \|_{\mathbb{L}^2(\mu)} \leq (b_i-a_i)
 \| L P_{a_i} f \|_{\mathbb{L}^2(\mu)}.
$$ 
Strong absolute continuity follows as, by spectral theory, for any $\alpha >0$ and any 
$f \in \mathbb{L}^2(\mu)$,
\begin{multline}
 \label{eq:borne_derivee_alpha}
\| L P_{\alpha} f \|_{\mathbb{L}^2(\mu)}^2 = \mu \Bigl( \bigl((-L) P_\alpha f  \bigr)^2\Bigr)\\
=\mu \Bigl( \bigl((-L)^2 P_{2\alpha} f\bigr) \, f \Bigr)
= \frac{1}{\alpha^2} \int_0^\infty (\alpha \xi )^2 e^{-2 \alpha \xi} \nu_{f}(d\xi) 
\leq  \frac{C}{\alpha^2} \mu( f^2),
\end{multline}
for some constant  $C>0$. 
Note that,  one also has 
$\mathcal{E}(P_\varepsilon f) \leq \frac{C}{\varepsilon} \mu(f^2)$ for any 
$f \in \mathbb{L}^2(\mu)$. It follows that 
$\| P_\varepsilon f \|_{\mathcal{D}} \leq \sqrt{1 + \frac{C}{\varepsilon}} 
\|f\|_{\mathbb{L}^2(\mu)}$. By duality, $\|f\|_{\mathbb{L}^2(\mu)}\leq \|f\|_{\mathcal{D}'}$, then 
$$
\| P_\varepsilon z \|_{\mathbb{L}^2(\mu)} \leq \sqrt{1 + \frac{C}{\varepsilon}} 
\|z\|_{\mathcal{D}'} \in \mathbb{L}^2([0,T]).
$$
We'll write $\tilde{z} \equiv P_\varepsilon z \in 
\mathbb{L}^2([0,T],\mathbb{L}^2(\mu))$ (or even sometime $\tilde{z} \equiv P_{\varepsilon/2} z$).

\noindent
We now turn our attention to the second term, 
$$
\Psi_\varepsilon(z)(t) \equiv \int_0^t P_{t-s+\varepsilon} z(s) ds, \quad (\varepsilon > 0). 
$$
First,  we show absolute continuity on $[0,T]$ of $\Psi_\varepsilon(z)$
in $\mathbb{L}^2(\mu)$. With $([a_i,b_i])_{i=1, \dots, N}$  a finite collection of non overlapping 
subintervals of $[0,T]$, 
\begin{multline*}
\| \Psi_\varepsilon(z)(b_i) - \Psi_\varepsilon(z)(a_i) \|_{\mathbb{L}^2(\mu)} =\\
 \| \int_{a_i}^{b_i}  P_{b_i-s+ \varepsilon}(z(s)) ds  + \int_0^{a_i} ds 
[P_{b_i-s} - P_{a_i-s}](P_\varepsilon z(s)) \|_{\mathbb{L}^2(\mu)} \\
\leq \int_{a_i}^{b_i} \| P_\varepsilon z(s) \|_{\mathbb{L}^2(\mu)} ds 
 +  
\int_0^{a_i} ds \int_{a_i-s}^{b_i-s} \|P_{\tau} L  (P_\varepsilon z(s)) \|_{\mathbb{L}^2(\mu)}\\
\leq \int_{a_i}^{b_i} \| \tilde{z}(s) \|_{\mathbb{L}^2(\mu)} ds 
 +  
\frac{C}{\varepsilon} (b_i -a_i) \int_0^{T}   
\| z(s) \|_{\mathcal{D}'} ds
\end{multline*}
by another use of \eqref{eq:borne_derivee_alpha}. (Strong) Absolute continuity follows.

\medskip
Continuity of $u$ at $t=0$ in $\mathbb{L}^2(\mu)$ follows by $C^0$ property of
the semigroup. Indeed, 
$\| \Psi_\varepsilon(z)(t) \|_{\mathbb{L}^2(\mu)}   
\leq \int_0^t \|  \tilde{z}(s) \|_{\mathbb{L}^2(\mu)} ds$
which goes to $0$ as $t$ goes to $0$. 

\medskip
To prove that, for
any $t$ a.e.,  $ \Psi_\varepsilon(z)(t) \in \mathcal{D}$ and  
$\int_0^T \mathcal{E}(\Psi_\varepsilon(z)(t)) dt < +\infty$, we will use spectral theory.
One has
\begin{multline*}
 \mathcal{E}(\int_0^t P_{t-s+\varepsilon}(z(s)) ds, \int_0^t P_{t-\tau+\varepsilon}(z(\tau)) d\tau ) 
\leq  \int_{[0,t]^2} \! \!  \! \! \! \!  \! \! ds \, d\tau \, 
\mathcal{E}( P_{t-s+\varepsilon}(z(s)),  P_{t-\tau+\varepsilon}(z(\tau)))\\
\leq t \int_0^t \mathcal{E}( P_{t-s+\varepsilon}(z(s))) ds 
\end{multline*}
thanks to Cauchy-Schwarz inequality. 
Hence,
$$
\int_0^T \mathcal{E}( \Psi_\varepsilon(z)(t)) dt \leq T \, \int_0^T dt \int_0^t ds \, 
\mathcal{E}( P_{t-s}(\tilde{z}(s))) = T \, \int_0^T ds \int_s^T dt \, \mathcal{E}( P_{t-s}(\tilde{z}(s)))
$$
by a first use of Fubini-Tonelli theorem. Then, noting that 
$$
\mathcal{E}( P_{t-s}(\tilde{z}(s))) = \int_0^\infty \xi e^{-2(t-s) \xi} \, 
\nu_{\tilde{z}(s)}(d\xi)
$$
and using once again Fubini-Tonelli theorem for $dt \otimes \nu_{\tilde{z}(s)}(d\xi)$ at 
fixed $s$, one gets
\begin{multline*}
\int_0^T \mathcal{E}( \Psi_\varepsilon(z)(t)) dt \leq  T \, \int_0^T ds 
\int_0^\infty \nu_{\tilde{z}(s)}(d\xi) 
\, \xi \int_s^T e^{-2(t-s) \xi} dt \\ 
\leq T \, \int_0^T ds \int_0^\infty \nu_{\tilde{z}(s)}(d\xi) 
\, \xi \frac{1- e^{-{2 (T-s) \xi}}}{2 \xi} \leq \frac{T}{2} \int_0^T ds \mu(\tilde{z}^2(s)) < \infty.
\end{multline*}

\bigskip
Now, we show that provided $\varepsilon >0$, $\Psi_\varepsilon(z)(t)$
is differentiable in $\mathbb{L}^2(\mu)$ for any $t$ a.e. in $[0,T]$ and,  
$$
\forall t \text{ a.e, }  \frac{\partial}{\partial t} \Psi_\varepsilon(z)(t) = P_\varepsilon (z(t)) + \int_0^t P_{t-s} L P_\varepsilon z(s) ds = 
P_\varepsilon (z(t)) + L ( \Psi_\varepsilon(z)(t)). 
$$
Let $h>0$ (the case when $h<0$ 
can be dealt with in the same way). 
Let us consider (in $\mathbb{L}^2(\mu)$) the difference between the associated differential ratio
and the expected derivative
$$
\frac{1}{h} \Bigl[ \int_0^{t+h} P_{t+h-s+\varepsilon} z(s) ds -  \int_0^{t} P_{t-s+\varepsilon} z(s) ds 
\Bigr] - \int_0^{t} P_{t-s} L P_{\varepsilon} z(s) ds - P_\varepsilon z(t).
$$
We split it into three terms. 
First,
$$
(I) = \frac{1}{h}  \int_t^{t+h} \Bigl[ P_{t+h-s} - Id 
\Bigr] P_{\varepsilon} z(s) ds. 
$$
Secondly,
$$
(II)=   \frac{1}{h}  \int_t^{t+h} P_{\varepsilon} z(s) ds - P_\varepsilon z(t).
$$
And third,
$$
(III) = \int_0^t ds \Bigl[ \frac{(P_{t+h-s+\varepsilon}-P_{t-s+\varepsilon})}{h} \bigl(z(s)\bigr) - 
P_{t-s} L P_{\varepsilon} z(s) \Bigr].
$$
Now, these three terms all go to $0$ in $\mathbb{L}^2(\mu)$ as $0<h$ goes to $0$.

\bigskip
Indeed, we deal with the first term as for absolute continuity of 
$
\Psi_\varepsilon ( z)
$
above.
One has
\begin{multline*}
\Bigl\| \frac{1}{h} \int_t^{t+h}
 \left[ P_{t-s + h} - Id \right] 
\, P_{\varepsilon} z(s) ds \Bigr\|_{\mathbb{L}^2(\mu)} 
\leq \frac{1}{h} \int_t^{t+h} ds \int^{t+h-s}_{0}
\Bigl\| P_\tau L P_{\varepsilon} z(s) \Bigr\|_{\mathbb{L}^2(\mu)} d\tau\\
\leq \frac{C}{\varepsilon} \int_t^{t+h}   \frac{t+h-s}{h} \|\tilde{z}(s)\|_{\mathbb{L}^2(\mu)} ds 
\leq \frac{C}{\varepsilon} \int_t^{t+h} ds \|\tilde{z}(s)\|_{\mathbb{L}^2(\mu)},
\end{multline*}
which goes to 0 as $h \rightarrow 0$.

\bigskip
Convergence of $(II)$ to $0$ in $\mathbb{L}^2(\mu)$, and this for any $t$ a.e., follows from the easy part of the 
fundamental theorem of calculus for Bochner integrable functions with values in $\mathbb{L}^2(\mu)$ 
(proved via comparison with strongly Henstock-Kurzweil integrable functions and Vitali
covering arguments in \cite[Theorems 7.4.2 and 5.1.4.]{schwabik-ye}
for instance). 

%

\bigskip
Finally, we focus  on $(III)$.  
For any $s$ a.e., as $0<h$ goes to $0$,
$$
\frac{(P_{h+\varepsilon}-P_{\varepsilon})}{h} \bigl(z(s)\bigr) \rightarrow L P_{\varepsilon} z(s)
$$
in $\mathbb{L}^2(\mu)$ as $P_{\varepsilon} z(s) \in \mathcal{D}(L)$. And we can use dominated
convergence theorem as, for $g_\varepsilon(\tau, s) \equiv P_\tau ( P_\varepsilon z(s) )$,
$$
 \| \frac{\partial}{\partial  \tau} g_\varepsilon(\tau, s) \|_{\mathbb{L}^2(\mu)}^2 
= \| P_\tau (-L) P_\varepsilon z(s) \|_{\mathbb{L}^2(\mu)}^2 
\leq \frac{C}{\varepsilon^2} \| \tilde{z}(s) \|_{\mathbb{L}^2(\mu)}^2
$$
still using \eqref{eq:borne_derivee_alpha}.

\bigskip
At the end of the day, $u$ is a solution a.e. of \eqref{eq:equation_generale_decorelee}.
Deducing that $u$ is a weak solution is easy. If 
$\phi \in C^\infty([0,T], \mathbb{L}^2(\mu))$, by bilinearity, $u \phi$ is 
absolutely continuous in  $\mathbb{L}^1(\mu)$ on $[a,T]$, $0<a<T$, and so is the 
real valued function
$t \mapsto \mu (u(t) \phi(t))$. 
The weak formulation follows when $a \rightarrow 0$ in the integration by parts formula
$$
\int_a^tds\mu(\partial_s u \, \phi(s))=\mu(u(t)\phi(t))-\mu(u(a)\phi(a))-\int_a^tds\mu(u(s) \, 
\partial_s\phi).
$$
The proof is complete.
\end{eproof}

\subsection{A mollified problem}

\label{section:mollified_problem}

\begin{remark}
\label{remark:definition_T}
 In sections \ref{section:mollified_problem} to \ref{section:removing_smoothing} below, we use  notation introduced in the statement
of lemma \ref{lem-lineaire}. So $T>0$ is fixed, 
$A(t) \in \mathbb{L}^\infty([0,T], \mathbb{L}^{\Phi_2}(\mu))$
and  
$B(t) \in 
\mathbb{L}^2([0,T], \mathbb{L}^2(\mu))$.
\end{remark}

Let us fix $\varepsilon >0$ and let us consider 
the following mollified problem 
\begin{equation}
\label{eq:mollified_linear_pb}
\tag*{(\text{$\mathbb{CS_\varepsilon}$})}
 \left\{ \begin{array}{l}
          \partial_t u^{(\varepsilon)}(t) = L u^{(\varepsilon)}(t)+ P_\varepsilon \Bigl( - A(t) \, u^{(\varepsilon)}(t) 
+  B(t) \Bigr),\\
           u^{(\varepsilon)}(0)= f, f \in \mathbb{L}^{2}(\mu)
         \end{array}
\right.
\end{equation}

We will prove that, for any $\varepsilon>0$ (and 
with some more work still at the limit $\varepsilon\to 0$), 
the problem~{\it \ref{eq:mollified_linear_pb}} has a weak solution 
in $[0,T]$ that is 
${u^{(\varepsilon)}} \in \mathbb{L}^2([0,T], \mathcal{D})\cap C([0,T],\mathbb L^2(\mu))$ 
 and, for any $\phi \in {C}^\infty([0,T], \mathcal{D})$, and any $0 \leq t \leq T$, 
\begin{multline}
\label{eq:weak_WMP}
\tag*{(\text{weak-$\mathbb{CS_\varepsilon}$})}
- \int_0^t  \mu ( u^{(\varepsilon)}(s) \partial_s \phi (s) ) \, ds  +    
\mu \left( u^{(\varepsilon)}(t) \phi(t)  - u^{(\varepsilon)}(0) \phi(0)\right)
= - \int_0^t    \mathcal{E}(u^{(\varepsilon)}(s), \phi(s) ) \, ds \\
+  \int_0^t  \mu \left(\phi(s) P_{\varepsilon}\big[-A(s){u^{(\varepsilon)}(s)}+B(s)\big] \right) \, ds. 
\end{multline}

\bigskip
To handle this problem, let us consider the following iteration
scheme which, as we will prove later, converge to the unique weak solution 
$u^{(\varepsilon)}$ of our problem (\text{$\mathbb{CS_\varepsilon}$}).
Initially,
$$
\left\{
\begin{array}{rcl}
\partial_t u^{(\varepsilon)}_{0}
& = & Lu^{(\varepsilon)}_{0}   \\
{u^{(\varepsilon)}_0}_{|t=0} & = & f
\end{array}
\right.
$$
and then define 
\begin{equation} \label{eq:aprxsln}
u^{(\varepsilon)}_{n+1}(t) \equiv P_{t}f + \int_0^t  
 P_{\varepsilon + t-s} (-A(s) u^{(\varepsilon)}_{n}(s) + B(s))ds.
\end{equation}
It follows from Lemmas ~\ref{lemma:regularity_nonlinearity} and  \ref{lemma_plus-hypothetic} 
that, for any $f \in \mathbb{L}^2(\mu)$,
 $u_{n+1} ^{(\varepsilon)}\in C([0,T], \mathbb{L}^2(\mu) \cap \mathbb{L}^2([0,T], \mathcal{D})$, and that 
for any $t$ a.e. in $[0,T]$, $u_{n+1}^{(\varepsilon)}(t)$ is differentiable in $ \mathbb{L}^2(\mu)$ and 
\begin{equation}
\label{eq:aprox} 
\left\{
\begin{array}{rcl}
\partial_t u^{(\varepsilon)}_{n+1}
& = & L u^{(\varepsilon)}_{n+1}(t) + P_\varepsilon \Bigl( -
A(t) u^{(\varepsilon)}_{n}(t) + B(t) \Bigr),\\
{u^{(\varepsilon)}_{n+1}}_{|t=0} & = & f,
\end{array}
\right. 
\end{equation}


The convergence scheme we detail below is adapted from the one presented 
in \cite{FRZ} in another context.

\begin{proposition}[Uniform bound]
 \label{Prop.Unifbnd.First.Shell}
Fix $\varepsilon >0$ and $f \in \mathbb{L}^2(\mu)$. 
Let $u_n^{(\varepsilon)}$ be the recursive solution of
the mollified problem introduced above.

There exists $\beta \in (0, +\infty)$ and $0<T_0 \leq T$ both independent of $\varepsilon$ 
and of the initial condition $f$ such that  for any
$n\in\mathbb{N}$,
\begin{equation}
\label{eq:uniform_bound} 
\sup_{0\leq t\leq T_0}\left( \mu ((u_n^{(\varepsilon)})^2(t))
+ \int_0^t\mathcal E(u_n^{(\varepsilon)})(s) ds\right) \leq
\beta \big(\mu(f^2)+
||B(\cdot)||^2_{\mathbb{L}^2([0.T],\mathbb{L}^2(\mu))}\big).
\end{equation}
\end{proposition}

\begin{eproof}
We use the notation $\tilde{u}^{(\varepsilon)}_{n}=P_\varepsilon u^{(\varepsilon)}_{n}$. 
For any $t$ a.e., 
\begin{multline*}
 \frac{1}{2} \frac{d}{dt} \mu((u^{(\varepsilon)}_{n+1})^2)=
\mu(u^{(\varepsilon)}_{n+1}Lu^{(\varepsilon)}_{n+1})
-\mu(A(t)\tilde{u}^{(\varepsilon)}_{n+1}u^{(\varepsilon)}_{n})
+\mu(B(t)\tilde{u}^{(\varepsilon)}_{n+1})\\
 \leq -\mathcal E(u^{(\varepsilon)}_{n+1})+
\frac{1}{2}
\mu(|A(t)|((\tilde{u}^{(\varepsilon)}_{n+1})^2+(u^{(\varepsilon)}_{n})^2)
+\bigl(\mu(B^2(t))\bigr)^{\frac{1}{2}} \, 
\bigl(\mu(\tilde{u}^{(\varepsilon)}_{n+1})^2\bigr)^{\frac{1}{2}} \\
\leq -\mathcal E(u^{(\varepsilon)}_{n+1})+
\frac{1}{2}
\mu(|A(t)|((\tilde{u}^{(\varepsilon)}_{n+1})^2+(u^{(\varepsilon)}_{n})^2)
+\frac{1}{2}\bigl( \frac{1}{\gamma} \mu(\tilde{u}^{(\varepsilon)}_{n+1})^2 + 
\gamma \mu(B^2(t))\bigr). 
\end{multline*}
Let 
\begin{equation}
\label{def-mgamma}
M_\gamma \equiv ||\mu(e^{\gamma |A|(\cdot)})||_{\mathbb L^\infty([0,T])}.
\end{equation}
Note that 
$1 \leq M_\gamma<\infty$ for any $\gamma>0$ 
since $A \in \mathbb{L}^\infty([0,T], \mathbb{L}^{\Phi_2}(\mu))$.

By a similar argument, the entropic  and the logarithmic Sobolev inequalities give 
\begin{multline*}
\mu(|A(t)|(\tilde{u}^{(\varepsilon)}_{n+1})^2)\leq \frac1\gamma 
\ent_\mu((\tilde{u}^{(\varepsilon)}_{n+1})^2)
+\frac{\mu((\tilde{u}^{(\varepsilon)}_{n+1})^2)}{\gamma}\log\mu(e^{\gamma |A(t)|})\\
\leq
\frac{C_{LS}}{\gamma} \mathcal E(\tilde{u}^{(\varepsilon)}_{n+1})+
\frac{\mu((\tilde{u}^{(\varepsilon)}_{n+1})^2)}{\gamma}\log M_\gamma,
\end{multline*}
and similarly for the other term. So that
\begin{multline*}
\frac{1}{2} \frac{d}{dt} \mu((u^{(\varepsilon)}_{n+1}))^2 \leq 
-\mathcal E(u^{(\varepsilon)}_{n+1})+\frac{ C_{LS}}{2\gamma}[\mathcal 
E(\tilde{u}^{(\varepsilon)}_{n+1})+\mathcal E({u}^{(\varepsilon)}_{n})]\\
+\frac{1 + \log M_\gamma }{2\gamma}[\mu((\tilde{u}^{(\varepsilon)}_{n+1})^2)
+\mu(({u}^{(\varepsilon)}_{n})^2)] + 
\frac{\gamma}{2} \mu(B^2(t)).
\end{multline*}
Using  $\mathcal E(\tilde{u}^{(\varepsilon)}_{n+1})
\leq\mathcal E({u}^{(\varepsilon)}_{n+1})$ 
and 
$\mu((\tilde{u}^{(\varepsilon)}_{n+1})^2)
\leq \mu(({u}^{(\varepsilon)}_{n+1})^2)$ 
and integrating with respect to $t$, 
\begin{multline*}
\mu(({u}^{(\varepsilon)}_{n+1})^2(t))+2(1-\frac{ C_{LS}}{2\gamma})
\int_0^t\mathcal E({u}^{(\varepsilon)}_{n+1})(s)ds
\leq \mu(f^2)\\
+\frac{1+ \log M_\gamma}{\gamma}
\int_0^t\mu(({u}^{(\varepsilon)}_{n+1})^2)(s)ds+
\frac{1+\log M_\gamma}{\gamma}
\int_0^t\mu(({u}^{(\varepsilon)}_{n})^2)(s)ds\\
+\frac{C_{LS}}{\gamma}
\int_0^t \mathcal{E}({u}^{(\varepsilon)}_{n})(s)ds+
\gamma ||B(\cdot)||^2_{\mathbb{L}^2([0.T],\mathbb{L}^2(\mu))}.
\end{multline*}
Choosing $\gamma> \frac{C_{LS}}{2}$, 
$\kappa_\gamma \equiv 1-\frac{ C_{LS}}{2\gamma}>0$ 
and setting
$$
\theta_{n}(t)=\mu(({u}^{(\varepsilon)}_{n})^2)(t)
+2\kappa_\gamma \int_0^t\mathcal E({u}^{(\varepsilon)}_{n})(s)ds,
$$
the above inequality implies 
\begin{multline*}
\theta_{n+1}(t)\leq \mu(f^2)+\gamma ||B(\cdot)||^2_{\mathbb L^2([0.T],\mathbb{L}^2(\mu))}
+\frac{1+\log M_\gamma}{\gamma}\int_0^t\theta_{n+1}(s)ds\\+
\frac{1+\log M_\gamma}{\gamma}\int_0^t\theta_{n}(s)ds+
\frac{ C_{LS}}{2\gamma\kappa_\gamma} \theta_{n}(t).
\end{multline*}
Hence, by Gronwall type arguments, one gets 
$$
\theta_{n+1}(t)
\leq e^{\frac{1+\log M_\gamma}{\gamma}t}
\Big[\alpha+\frac{1+\log M_\gamma}{\gamma}
\int_0^t\theta_n(s)ds+
\frac{ C_{LS}}{2\gamma\kappa_\gamma}\theta_n(t)\Big],
$$
where 
$$
\alpha=
\mu(f^2)+\gamma||B(\cdot)||^2_{\mathbb L^2([0.T],\mathbb{L}^2(\mu))}.
$$
It gives, for any $0<T_0 \leq T$,
$$
\sup_{t\in[0,T_0]}\theta_{n+1}(t)\leq e^{\frac{1+\log M_\gamma}{\gamma}T_0} 
\alpha 
+e^{\frac{1+\log M_\gamma}{\gamma}T_0}
\Big[\frac{ 1+ \log M_\gamma}{\gamma}T_0
+\frac{ C_{LS}}{2\gamma - C_{LS}}\Big]
\sup_{t\in[0,T_0]}\theta_{n}(t)
$$
Let us denote $Z_n=\sup_{t\in[0,T_0]}\theta_{n}(t)$. 

Now, provided we choose $\gamma >  C_{LS}$, 
$\frac{ C_{LS}}{2 \gamma - C_{LS}}<1$,
so that, for $T_0>0$ small enough, 
$$
\eta_{T_0}=e^{\frac{1+\log M_\gamma}{\gamma}T_0}
\Big[\frac{1 +\log M_\gamma}{\gamma}T_0
+\frac{C_{LS}}{2\gamma - C_{LS}}\Big]<1. 
$$
we end up with 
$$
Z_{n+1}\leq e^{\frac{1+\log M_\gamma}{\gamma}T_0} 
\alpha + \eta_{T_0} Z_n,
$$
Hence, by induction, 
$$
Z_n\leq \alpha e^{\frac{1+\log M_\gamma}{\gamma}T_0} 
(1+\cdots+\eta_{T_0}^{n-1})+ \eta_{T_0}^n
Z_0.
$$
Note that 
\begin{multline*}
Z_0=\sup_{t\in[0,T_0]}\big\{\mu(P_{t}(f)^2)
+2\kappa_\gamma \int_0^t\mathcal E(P_{s}(f)) ds\big\}\\
\leq \sup_{t\in[0,T_0]}\big\{\mu(P_{t}(f)^2)
+2\int_0^t\mathcal E(P_{s}(f))ds\big\}\leq \mu(f^2)\leq \alpha,
\end{multline*}
since the map $s\mapsto \mu(P_{t}(f)^2)
+2\int_0^t\mathcal E(P_{s}(f))ds$ is decreasing. 
It follows that, for any $n\geq0$, 
$$
Z_n\leq \alpha e^{\frac{1+\log M_\gamma}{\gamma}T_0} (1+\cdots+\eta_{T_0}^{n})
\leq \alpha e^{\frac{1+\log M_\gamma}{\gamma}T_0} \frac{1}{1-\eta_{T_0}}, 
$$
which is the expected bound.  
\end{eproof}

\begin{proposition}[Existence for mollified problem; $\varepsilon>0$] \label{Prop.Convergence.First.Shell}
For any $\varepsilon>0$ and any initial datum $f \in \mathbb{L}^2(\mu)$, 
there exists a weak solution  
$u^{(\varepsilon)}$ on $[0,T]$  of the mollified 
problem~{\it \ref{eq:mollified_linear_pb}} as defined in~\ref{eq:weak_WMP}.
\end{proposition}

\begin{eproof}
Let $w_{n+1}^{(\varepsilon)}=u_{n+1}^{(\varepsilon)}-u_{n}^{(\varepsilon)}$ and 
$\tilde{w}_{n}^{(\varepsilon)}=P_\varepsilon (w_{n}^{(\varepsilon)})$.
For any $t\geq0$ a.e., 
$$
\frac12\frac{d}{dt}\mu((w_{n+1}^{(\varepsilon)})^2)=-\mathcal E(w_{n+1}^{(\varepsilon)})-\mu(\tilde{w}_{n+1}^{(\varepsilon)}A(t)w_n^{(\varepsilon)}).
$$
Again thanks to the entropic and the logarithmic Sobolev inequalities, 
\begin{multline*}
\frac12\frac{d}{dt}\mu((w_{n+1}^{(\varepsilon)})^2)\leq-\mathcal E(w_{n+1}^{(\varepsilon)})+\frac{ C_{LS}}{2\gamma}\big[\mathcal E(\tilde{w}_{n+1}^{(\varepsilon)})+\mathcal E({w}_{n}^{(\varepsilon)})\big]\\
+\frac{\log M_\gamma}{2\gamma}\big[\mu((\tilde{w}_{n+1}^{(\varepsilon)})^2)+\mu(({w}_{n}^{(\varepsilon)})^2)\big],
\end{multline*}
where $M_\gamma$ were defined in the proof of Proposition~\ref{Prop.Unifbnd.First.Shell}. By the same arguments as before,
\begin{multline*}
\mu((w_{n+1}^{(\varepsilon)})^2(t))+2\kappa_\gamma\int_0^t\mathcal E(w_{n+1}^{(\varepsilon)})(s)ds\leq 
\frac{\log M_\gamma}{\gamma}\int_0^t\mu((w_{n+1}^{(\varepsilon)})^2)(s)ds\\
+\frac{\log M_\gamma}{\gamma}\int_0^t\mu((w_{n}^{(\varepsilon)})^2)(s)ds+\frac{C_{LS}}{\gamma}\int_0^t\mathcal E(w_{n}^{(\varepsilon)})(s)ds,
\end{multline*}
 with again $\kappa_\gamma=1-\frac{C_{LS}}{2\gamma}>0$ provided we choose $\gamma >  
\frac{C_{LS}}{2}$. 

Fixing $0<\tilde{T}_0\leq T_0$, where $T_0$ has been 
defined in the previous proposition,
and mimicking what we have done to prove that proposition, this leads to 
$$
\sup_{t\in[0,\tilde{T}_0]}\{\tilde{\theta}_{n+1}(t)\}
\leq \tilde{\eta}_{\tilde{T_0}}
 \sup_{t\in[0,T]}\{\tilde{\theta}_{n}(t)\},
$$
where 
 $
\tilde{\theta}_{n}(t)=\mu((w_{n+1}^{(\varepsilon)})^2)
+2\kappa_\gamma\int_0^t\mathcal E(w_{n+1}^{(\varepsilon)})(s)ds
$ 
and where 
$$
\tilde{\eta}_{\tilde{T}_0}= (M_\gamma)^{\frac{\tilde{T}_0}{\gamma}}
\Big[\frac{\log M_\gamma}{\gamma} \tilde{T}_0
+\frac{C_{LS}}{2\gamma - C_{LS}}\Big]. 
$$
If we choose $\gamma > C_{LS}$,
we may take $0<\tilde{T}_0\leq T_0$ small enough (and independent of the initial 
condition $f$)
so that $\tilde{\eta}_{\tilde{T}_0}<1$.

Iterating and using uniform bound \eqref{eq:uniform_bound} for $n=1$ (and $n=0$), one
gets 
$$
\sup_{t\in[0,\tilde{T}_0]}\{\tilde{\theta}_{n+1}(t)\} \leq \tilde{\beta} 
\big(\mu(f^2)+||B(\cdot)||^2_{\mathbb L^2([0.T],\mathbb{L}^2(\mu))}\big) 
\, \tilde{\eta}_{\tilde{T}_0}^n.
$$
It follows that
$(u_n^{(\varepsilon)})_{n\in\N}$ 
is  a Cauchy sequence in $\mathbb L^2([0,\tilde{T}_0],\mathcal D)\cap 
C([0,\tilde{T}_0],\mathbb L^2(\mu))$. 

It converges to some $u^{(\varepsilon)}$ which is a weak solution in $[0,\tilde{T}_0]$ 
of~{  \ref{eq:mollified_linear_pb}} (see 
page~\pageref{weak-solution}, but note that things are much simpler here). As $\tilde{T}_0$
does not depend on $f$, one easily extends the solution to the entire interval $[0,T]$.   
\end{eproof}

\subsection{Uniqueness}
We now state uniqueness of a weak solution for both cases : 
with or without a mollification.

\begin{proposition}[Uniqueness]
For any $\varepsilon \geq 0$, a weak solution  
$u^{(\varepsilon)}$ on $[0,T]$  of the 
problem~{$\ref{eq:mollified_linear_pb}$} with initial datum $f \in \mathbb{L}^2(\mu)$ is unique. 
\end{proposition}

We omit the proof which is quite similar to the one of proposition
\ref{prop-unique}.

\subsection{Existence for the cornerstone linear problem}

\label{section:removing_smoothing}

Recall remark \ref{remark:definition_T}: \\
$T>0$ is fixed and $A(t) \in \mathbb{L}^\infty([0,T], \mathbb{L}^{\Phi_2}(\mu))$
and  
$B(t) \in  \mathbb{L}^2([0,T], \mathbb{L}^2(\mu))$.

\begin{proposition}[Removing the smoothing]
Let $f \in \mathbb{L}^2(\mu)$.
There exists $0<T_0 \leq T$ (independent of $f$) such that the weak 
solution $u^{(\varepsilon)}$, $\varepsilon >0$, of the mollified problem~{$\ref{eq:mollified_linear_pb}$}, 
(with the same initial datum $f$)
converges as $\varepsilon$ goes to 0, to some limit function 
$u$ in  $\mathbb{L}^2([0,T_0], \mathcal{D})\cap C([0,T_0],\mathbb L^2(\mu))$.
Moreover, $u$ may be extended to a weak solution of the cornerstone linear  
problem~{\it \ref{eq:cornerstone_linear}}, with initial datum $f$, on $[0,T]$.
\end{proposition}

\begin{eproof}
Let $\varepsilon_1>\varepsilon_0>0$ and let $u_0=u^{(\varepsilon_0)}$ 
and $u_1=u^{(\varepsilon_1)}$
be the associated solutions of the mollified problem~\ref{eq:weak_WMP}. 
Using Steklov calculus as in the previous proof, 
we get the same estimate as if we were dealing with strong solutions. Here we avoid 
such technicalities to focus on the main arguments.   
Let us denote $w=u_1-u_0$ and $\tilde{w} = P_{\varepsilon_1} w$. 
One has
\begin{multline*}
\frac12\frac{d}{dt}\mu(w^2)=-\mathcal E(w)+\mu(w[P_{\varepsilon_1}(-A(t)u_1+B(t))-P_{\varepsilon_0}(-A(t)u_0+B(t))])\\
=-\mathcal E(w)
+ \underbrace{\mu \Bigl(w P_{\varepsilon_1}(-A(t)w) \Bigr)}_{(I)} + 
\mu \Bigl(w(P_{\varepsilon_1}-P_{\varepsilon_0})(-A(t)u_0+B(t))\Bigr).
\end{multline*}
Term $(I)$ is bounded by $\frac{C_{LS}}{\gamma} \mathcal{E}(w(t))
+ \frac{\log(M_\gamma)}{\gamma} \mu (w^2(t))$ as in the previous proof.
After integration, using symmetry of the semigroup, one gets
\begin{multline*}
\mu(w^2(t)) + (1 - \frac{C_{LS}}{\gamma}) \, \int_0^t \mathcal{E}(w(s)) ds\\
\leq \frac{\log(M_\gamma)}{\gamma} \int_0^t \mu (w^2(s)) ds
+ \underbrace{\int_0^t \mu \Bigl((P_{\varepsilon_1}-P_{\varepsilon_0})(w) \, (-A(s)u_0(s)+B(s))\Bigr)ds}_{(II)},
\end{multline*}
(which is the estimate we would get rigorously after letting $h \rightarrow 0$ in the Steklov regularisation). 
After using Gronwall type arguments and taking the supremum over $t \in [0,T_0]$, $0<T_0 \leq T$, 
we note that, if we prove term $(II)$
goes to $0$ as $\varepsilon_1>\varepsilon_0>0$ both go to $0$, 
then $(u^{(\varepsilon)})_{\varepsilon>0}$
is Cauchy (as $\varepsilon$ goes to $0$)
in the Banach space 
$\mathbb{L}^2([0,T_0], \mathcal{D})\cap C([0,T_0],\mathbb L^2(\mu))$.
Now, by Cauchy-Schwarz inequality,
\begin{equation}
\label{eq:difference_bound}
(II) \leq \left(\int_0^t ds \mu \left[\left(
P_{\varepsilon_1}-
P_{\varepsilon_0}\right)(w(s))\right]^2\right)^\frac12
\left(\frac12 \int_0^t ds \mu \bigl(A^2(s) u_0^2(s) + B^2(s)\bigr)\right)^\frac12. 
\end{equation}
Following lemma \ref{lemma:regularity_nonlinearity},
\begin{multline*}
\int_0^t ds \mu \bigl(A^2(s) u_0^2(s) + B^2(s)\bigr)
\leq \| A u_0 \|_{\mathbb{L}^2([0,T], \mathbb{L}^2(\mu))}^2 + 
 \| B \|_{\mathbb{L}^2([0,T], \mathbb{L}^2(\mu))}^2 \\ 
\leq  \max(\log(2), C_{LS}) \| A \|_{\mathbb{L}^\infty([0,T], \mathbb{L}^{\Phi_2}(\mu))}^2 \, 
\| u^{(\varepsilon_0)} \|_{\mathbb{L}^2([0,T], \mathcal{D})}^2 +
\| B \|_{\mathbb{L}^2([0,T], \mathbb{L}^2(\mu))}^2.
\end{multline*}
Choosing $T_0$ as in Proposition \ref{Prop.Unifbnd.First.Shell},
one may pass to the limit $n \rightarrow \infty$ in the uniform bound
\eqref{eq:uniform_bound} to get that, for any $\varepsilon >0$, 
\begin{equation}
 \label{eq:norm_uniform_bound}
\| u^{(\varepsilon)} \|_{\mathbb{L}^2([0,T], \mathcal{D})}^2 \leq 
\beta \, (T_0+1) \, (\mu(f^2) + \|B\|_{\mathbb{L}^2([0,T], \mathbb{L}^2(\mu))}^2.
\end{equation}
So the second factor of  \eqref{eq:difference_bound}
is bounded uniformly in $\varepsilon_0$.
In order to prove convergence to $0$ of the other factor $\int_0^t ds \mu \left[\left(
P_{\varepsilon_1}-
P_{\varepsilon_0}\right)(w(s))\right]^2$
when $\varepsilon_1>\varepsilon_0>0$ both go to $0$, one makes use of spectral theory
and the above uniform bound \eqref{eq:norm_uniform_bound}. Details are given in 
\cite[Theorem 4.10]{FRZ}. 

Eventually, the limit $u$ of  $(u^{(\varepsilon)})_{\varepsilon>0}$ (as $\varepsilon$ goes to $0$)
in  
$\mathbb{L}^2([0,T_0], \mathcal{D})\cap C([0,T_0],\mathbb L^2(\mu))$ is a weak solution, which 
can be extended to a weak solution  on the entire interval $[0,T]$ as $T_0$
doesn't depend on the initial datum $f$. 
\end{eproof}

\subsection{Non-negativity}

We prove here that, provided $A$ and $B$ are nonnegative,  the weak solution $u$
of problem ~{$\ref{eq:cornerstone_linear}$}, with a nonnegative 
initial datum $f$, is  nonnegative.

Let us define $u_{-}= (-u)_+ = \max(-u, 0)$. Then, formally,
\begin{multline*}
 \frac{1}{2} \frac{d}{dt} \mu( (u_{-}(t))^2) = 
- \mu( u_{-}(t) \partial_t u ) = -\mu((u_-(t)Lu(t)) + \mu(u_{-}(t) A(t) u(t)) - \mu( u_-(t) B(t) ) \\
\leq - \mu \Bigl( A(t)\underbrace{(-u)_+(t) (-u)(t)}_{= ((-u)_+(t))^2} \Bigr)+\mathcal E(u_-(t),u(t))
\leq 0 
\end{multline*}
using positivity of $A(\cdot)$ and $B(\cdot)$,
and 
$\mathcal E(u_-(t),u(t))\leq -\mathcal E(u_-(t),u_-(t))\leq 0$. 

\bigskip
Rigorous arguments to get this are as follows. We  consider the Steklov average
$a_h(u)(t)$ and its negative part $a_h^-(u)(t) \equiv
\max(0,-a_h(u)(t)).$ 
Recall that, as $h$ goes to $0$, for any $t \in [0,T]$, 
$a_h^-(u)(t) \rightarrow u^-(t)$ in
$\mathbb{L}^2(\mu)$ and $a_h(u) \rightarrow u$ in
$\mathbb{L}^2([0,T], \mathcal{D})$. 
It follows that
$a_h^{-}(u) \rightarrow u^{-}$ in
$\mathbb{L}^2([0,T], \mathcal{D})$.
Namely, from any sequence going to $0$, extract a subsequence $(h_n)$
such that, for any $t$ a.e. in $[0,T]$, 
$a_{h_n}(u)(t) \rightarrow u(t)$ in $\mathcal{D}$. By continuity of contractions
\cite{ancona}, it follows  $a_{h_n}^{-}(u)(t) \rightarrow u^{-}(t)$ , in $\mathcal{D}$, $t$
a.e. and one may check easily that the sequence 
$(\| a_{h_n}^{-}(u)(t) -u^{-}(t)\|^2_{\mathcal{D}})_n$
is uniformly integrable in $\mathbb{L}^1([0,T])$. 

Moreover, in $W^{1,2}((0,T),
\mathbb{L}^2(\mu))$,
$$
\partial_s a_h^-(u) (s) = -
\partial_s a_h(u) (s) \chi_{\{a_h(u)(s) \leq 0 \}} = - \frac{1}{h} (u(s+h) -
u(s)) \chi_{\{a_h(u)(s) \leq 0\}}
$$
where $\chi$ denotes the indicator function. Hence, using the
definition of a weak solution (with the constant test function
$a_h^{-}(u)(s) \in \mathcal{D}$), we get
\begin{multline*}
\frac 12\mu \left( a_h^{-}(u)(t) \right)^2 = \frac 12\mu \left(
a_h^{-}(u)(0) \right)^2 +
\frac 12\int_0^t ds \, \partial_s \mu \left( a_h^{-}(u)(s) \right)^2 \\
= \frac 12\mu \left( a_h^{-}(u)(0) \right)^2 - \int_0^tds \, \mu
\left( a_h^{-}(u)(s)\frac1h\left(u(s+h)-u(s)\right) \right)
\\
 = \frac 12\mu \left( a_h^{-}(u)(0) \right)^2  + \int_0^t  \! \! ds \frac1h
\int_s^{s+h}  \! \!   \! \! \! \! \! \! d\tau \left[ \mathcal{E}\left(a_h^{-}(u)(s),
u(\tau)\right)
 +   \mu\Bigl(a_h^{-}(u)(s) \, \bigl(A(\tau) u(\tau) -B(\tau) \bigr) \Bigr) \right]\\
= \frac 12 \mu \left( a_h^{-}(u)(0) \right)^2  + \int_0^t  \! \!  \! \! ds  \! \! \left[
\mathcal{E}\Bigl(a_h^{-}(u)(s), a_h(u)(s)\Bigr) +   
\mu\Bigl(a_h^{-}(u)(s) \, a_h\bigl(A(\cdot) u(\cdot) -B(\cdot) \bigr)(s) \Bigr)
 \right]
\end{multline*}
We can pass to the limit with $h\to 0$ 
which yields (as $\mu
\left((f^{-})^2 \right) =0$)
\begin{equation*}
\frac 12\mu \left( u^{-}(t) \right)^2  = \int_0^tds  \, \mathcal{E}(
u^{-}(s), u(s)) + \mu\Bigl(u^{-}(s) \, (A(s) u(s) -B(s) \bigr) \Bigr)
 \leq 0,
\end{equation*}  
for the same reason as above.

\medskip
The proof of Lemma~\ref{lem-lineaire} is complete.

\section{Extension to the general case}

\label{section:extension}

The chemical reactions we consider here are of the following form 
$$
\sum_{i\in F} \alpha_i\A_i \rightleftharpoons \sum_{i\in F} \beta_i\A_i,
$$
for some given integers $\alpha_i\neq\beta_i$, for any $i\in F$. 
$F=\{1, \dots, q \}$ is
a finite set. The associated reaction-diffusion
equation is (after appropriate change of variables)
\begin{equation}
\label{eq:RDP_extension}
\left\{
\begin{array}{l}
\partial_t u_i= L_i u_i + \lambda_i (\beta_i-\alpha_i) \,
\PAR{\prod_{j=1}^q u_j^{\alpha_j} -\prod_{j=1}^q u_j^{\beta_j}}\\
u_{i|t=0} = f_i, \quad i \in F
\end{array}
\right. 
\end{equation}
This equation is a particular form of the abstract equation \ref{eq:RDP} 
on page \pageref{eq:RDP} with constant vector 
$\lambda_i(\beta_i-\alpha_i)$, $i=1, \dots, q$ and nonlinearity
$G(\vec{u})=  \prod_{j=1}^q u_j^{\alpha_j} -\prod_{j=1}^q u_j^{\beta_j}$.
The method we detailed for the two-by-two case may be adapted to this 
general situation provided the following assumptions hold.

\subsubsection*{Linearity assumptions}

We assume that

\begin{enumerate}
 \item $F$ may be partitioned as $F = \sqcup_{k \in K} F_k$
so that,  $L_i$ only depends on which $F_k$, $k \in K$, the index $i$ 
belongs to. 
We denote by $\tilde{L}_k$ the common operator for any $i \in F_k$. 
For any $k \in K$, one has the following:
\item $\tilde{L}_k$ is a Markov generator with (selfadjoint in the $\mathbb{L}^2$
space associated with the) invariant probability measure $\mu_k$ on 
$(\mathbb{M}, \mathcal{B}_{\mathbb{M}})$ 
(with the same assumptions as in page
\pageref{section:abstract_problem_assumptions}).
\item $(\tilde{L}_k,\mu_k)$ satisfies logarithmic Sobolev inequality 
with constant $C_k$.
\item The measures $(\mu_k)_{k \in K}$ are mutually equivalent in the strong sense that 
there exists a measure $\mu$ on $(\mathbb{M}, \mathcal{B}_{\mathbb{M}})$ and 
$C \in (1, +\infty)$ such that
$$
\forall k \in K, \qquad \frac{1}{C} \leq \frac{d \mu_k}{d \mu} \leq C.
$$
\end{enumerate}

\subsubsection*{Nonlinearity assumptions} 
Let $F_{k}^-=\{i\in F_k,\,\,\beta_i-\alpha_i<0\}$ and $F_{k}^+=\{i\in F_k,\,
\,\beta_i-\alpha_i>0\}$. 
We assume that, for any $k \in K$, $F_{k}^-$ and $F_{k}^+$ are not empty.

\smallskip
(Note that this replaces, in the present context, the hypothesis
we made in the two-by-two case that $C_1=C_3$ and $C_2=C_4$.)

%

\subsubsection*{Initial data assumptions}

We assume the following common exponential
integrability on the initial data. 

\bigskip
\begin{minipage}{11cm}
\textbf{Common integrability assumption.} 
We assume  that, for any $i=1,\dots,q$,
$f_i \in E^{\Phi_{2 \theta}}(\mu)$, where
$\theta \equiv \max ( \sum_{i=1}^q \alpha_i ,  \sum_{i=1}^q \beta_i ) -1$.
\end{minipage}

\subsubsection*{Iterative sequence}

We now define an approximation sequence 
$(\vec{u}^{(n)}(t))_{n \in \mathbb{N}}$ which converges 
to the solution
of problem \eqref{eq:RDP_extension}. 
It is obtained recursively 
as solutions of the following linear problems.

\bigskip  
Let us fix a nonnegative 
initial datum $\vec{f}$ satisfying the integrability assumptions 
introduced before.  

\bigskip
For any $n \geq 0$, we will impose $\vec{u}^{(n)}(0) = \vec{f}$
and, for $n=0$, $\partial_t u_i^{(0)} = L_i u_i^{(0)}$, $i=1, \dots, q$.

\bigskip
\noindent
Let $N_k = |F_k|$, $N_k^+ = |F_k^+|$ and $N_k^- = |F_k^-|$.
Assume $N_k^- \geq N_k^+$ (the other case is similar by symmetry). Let us label elements of $F_k^\pm$ in the following way
$$
F_k^- = \{ i^{k-}_1, \dots , i^{k-}_{N_k^-} \} \qquad \text{ and } \qquad F_k^+ = \{ i^{k+}_1, \dots , 
i^{k+}_{N_k^+} \}.
$$
We consider an onto mapping $\nu_k : F_k^- \twoheadrightarrow F_k^+$ 
defined by
$$
\nu_k(i_l^{k-}) = i_m^{k+} \quad \text{ provided } l - m  \in N_k^+ \mathbb{Z}. 
$$
Define furthermore, 
for any $i,j \in F$, 
$
\alpha_j^{(i)}= \left\{
\begin{array}{ll}
\alpha_j & \text{if } j \neq i\\
\alpha_j-1 & \text{if } j = i
\end{array}
\right.
$
and similarly for $\beta$'s. Let us note here that, 
for any $i \in F_k^+$ and $j \in F_k^-$, 
$\beta_i >0$ and $\alpha_j >0$. 
Finally, let $\delta_i = \lambda_i |\beta_{i}-\alpha_{i}|>0$, for any $i \in F$.

\bigskip
\noindent
The iterated sequence is then defined as follows\footnote{We recommend to translate at first reading the following
general case in the simpler two-by-one case 
$\mathcal{A}_1 + \mathcal{A}_2 \rightleftharpoons \mathcal{A}_3$ with the same diffusion operator.}.
In the case $i \in F_k^-$, 
\begin{equation}
\label{eq:approx_sequence_general_->+}
\partial_t u_{i}^{(n)}= \tilde{L}_k u_{i}^{(n)} - \delta_i
\Bigl(
\prod_{j=1}^q (u_j^{(n-1)})^{\alpha_j^{(i)}} \, u_{i}^{(n)} -
\prod_{j=1}^q (u_j^{(n-1)})^{\beta_{j}^{(\nu_k(i))}} \, u_{\nu_k(i)}^{(n)} 
\Bigr).
\end{equation}
And, in the case $i \in F_k^+$,
\begin{equation}
\label{eq:approx_sequence_general_+<-} 
\partial_t u_{i}^{(n)}= 
\tilde{L}_k u_{i}^{(n)} + 
\frac{\delta_i}{Z_{k,i}} \sum_{r \in \nu_k^{-1}(i)} \delta_r
\Bigl( \prod_{j=1}^q (u_j^{(n-1)})^{\alpha_j^{(r)}} \, u_{r}^{(n)} -
\prod_{j=1}^q (u_j^{(n-1)})^{\beta_{j}^{(i)}} \, u_{i}^{(n)} 
\Bigr).
\end{equation}
where $Z_{k,i} = \sum_{r \in \nu_k^{-1}(i)} \delta_r$.

\subsubsection*{Why the sequence is well defined.}

Recall the sequence starts with the heat semigroups associated
to the $L_i$'s 
$$
u_i^{(0)}(t) \equiv e^{t L_i} f_i, \qquad i \in F.
$$
%
\noindent It follows from appendix \ref{appendix:Bochner-measurability}
that under our assumptions on $\vec{f}$, $u_i^{(0)} \in \mathbb{L}^\infty([0,T],E^{\Phi_{2\theta}}(\mu))$.
We hence assume we have proved, $\vec{u}^{(n-1)}$ is well defined, for some
$n \geq 1$ and that 
\begin{equation}
\label{eq:regularity_u^(n)}
\vec{u}^{(n-1)} \in (\mathbb{L}^\infty([0,T],E^{\Phi_{2\theta}}(\mu)))^q,
\end{equation}
for any 
$T \in (0,+\infty)$. 
Lemma \ref{lemma:Young} below ensures that, for any $i \in F$, the 
mapping
$
(\su_1, \dots, \su_q) \in (E^{\Phi_{2 \theta}}(\mu))^q \mapsto 
\prod_{j=1}^q \su_j^{\alpha_j^{(i)}} \in E^{\Phi_2}(\mu)
$ 
is continuous, so that 
\begin{equation*}
\label{eq:prod_Ephi2}
 \prod_{j=1}^q (u_j^{(n-1)})^{\alpha_j^{(i)}} \in 
\mathbb{L}^\infty([0,T], E^{\Phi_2}(\mu))
\end{equation*}
(and similarly
for $\beta$'s).
 
\begin{lemma}
\label{lemma:Young}
Assume $p_1, \dots p_q \geq r >0$ such that 
$\frac{1}{p_1} + \dots + \frac{1}{p_q} =\frac{1}{r}$. 
Let $\Phi(x) = \exp(|x|)-1$ and recall $\Phi_\beta (x) = \Phi(|x|^\beta)$, for
any $\beta \geq 1$. 
Then, for any $\alpha>0$ such that $\alpha r \geq 1$, the $q$-linear mapping
$$
(u_1, \dots, u_q) \in \mathbb{L}^{\Phi_{\alpha p_1}}(\mu) \times \dots \times 
\mathbb{L}^{\Phi_{\alpha p_q}}(\mu) \mapsto u_1 \dots u_q \in  \mathbb{L}^{\Phi_{\alpha r}}(\mu)
$$
is continuous.
\end{lemma}

\begin{eproof}
Assume $u_i \neq 0$, for all $i=1, \dots, q$ and denote 
$\gamma_i \equiv \|u_i\|^{-1}_{\mathbb{L}^{\Phi_{\alpha p_i}}}$. Then one has
$$
\forall i = 1, \dots, q, \qquad \mu \Bigl( \prod_{j=1}^q u_j^{\alpha_j} -\prod_{j=1}^q u_j^{\beta_j}\exp \bigl( |\gamma_i u_i|^{\alpha p_i} \bigr) \Bigr)
\leq 2.
$$ 
The result will follow if we show that 
$$
\mu \Bigl( \exp \bigl( |\gamma_1 \dots \gamma_q \, u_1 \dots u_q|^{\alpha r} \bigr) \leq 2.
$$
Recall Young inequality: for any $a_1, \dots, a_q \geq 0$, 
$$
\frac{1}{r} a_1^r \dots a_q^r \leq \frac{a_1^{p_1}}{p_1} + \dots + \frac{a_q^{p_q}}{p_q}.
$$ 
Hence, using also H\"older inequality,
\begin{multline*}
 \mu\Bigl(e^{(\gamma_1 |u_1|)^{\alpha r} \dots (\gamma_q |u_q|)^{\alpha r})}\Bigr)
\leq 
\mu \Bigl(e^{\frac{r}{p_1} (|\gamma_1 u_1|^{\alpha p_1})} \dots 
e^{\frac{r}{p_q} |\gamma_q u_q|^{\alpha p_q})} \Bigr)\\
\leq \mu \Bigl(e^{|\gamma_1 u_1|^{\alpha p_1})} \Bigr)^{\frac{r}{p_1}}
\dots \mu \Bigl(e^{|\gamma_q u_q|^{\alpha p_q})} \Bigr)^{\frac{r}{p_q}} 
\leq 2^{\frac{r}{p_1} + \dots + \frac{r}{p_q}} = 2. 
\end{multline*}
\end{eproof} 

To prove recursively that the sequence $(\vec{u}^{(n)})_n$ is well defined, 
we have to split the 
cornerstone existence lemma into the following two lemmas.

\begin{elem}[Matrix cornerstone existence lemma]
\label{lem-lineaire-matriciel}
Let $(L,\mu)$ be a Markov generator satisfying  logarithmic Sobolev inequality
with constant $C_{LS} \in (0,\infty)$.  
Let $T>0$ and $A=A(t)$ be an $N \times N$ matrix with coefficients in 
$\mathbb{L}^\infty([0,T], \mathbb{L}^{\Phi_2}(\mu))$ and 
$\vec{B} \in (\mathbb{L}^2([0,T], \mathbb{L}^2(\mu)))^N$.
Then the Cauchy problems 
\begin{equation}
\label{eq:matrix_cornerstone}
\tag*{(\text{$\mathbb{MCS}$})}
 \left\{ \begin{array}{l}
          \partial_t \vec{u}(t) = L \vec{u}(t) + A(t) \, \vec{u}(t) + \vec{B}(t),\\
           \vec{u}(0)= \vec{f}, \vec{f} \in (\mathbb{L}^2(\mu))^N
         \end{array}
\right.
\end{equation}
and 
\begin{equation}
\label{eq:matrix_cornerstone_+}
\tag*{(\text{$\mathbb{MCS}_+$})}
 \left\{ \begin{array}{l}
          \partial_t \vec{u}(t) = L \vec{u}(t) + A(t) \, \vec{u}_+(t) + \vec{B}(t),\\
           \vec{u}(0)= \vec{f}, \vec{f} \in (\mathbb{L}^2(\mu))^N
         \end{array}
\right.,
\end{equation}
with $\vec{u}_+ = ((u_1)_+, \dots (u_q)_+)$, both  
have a unique weak solution on $[0,\infty)$ 
\end{elem}

\noindent Note that we use that $u \mapsto u_+$ is a contraction
so that it contracts both the $\mathbb{L}^2(\mu)$ norm and the Dirichlet
form $\mathcal{E}$.

\smallskip

\noindent In the system defined by \eqref{eq:approx_sequence_general_->+}
and \eqref{eq:approx_sequence_general_+<-}
only blocks made of some $i \in F_k^+$ and $j$'s in $\nu_k^{-1}(i)$ interact.
We now focus on these coordinates. 
The following lemma ensures that positivity and Bochner 
measurability \eqref{eq:regularity_u^(n)} propagate along 
the approximation
sequence.

\begin{elem}[Positivity and propagation of measurability.]
\label{lemma:approx_sequence_positivity_general}
Let $N \geq 2$ and let $\delta_1, \dots, \delta_{N-1} \geq 0$
such that $Z \equiv \sum_{i=1}^{N-1} \delta_i >0$. Assume 
furthermore $\vec{B}(t) = \vec{0}$ and $A(t)$ is of the following
form
\begin{equation}
\label{eq:specific_matrix}
A(t)= \left( \begin{array}{cccccc} 
- a_1(t) & 0 & 0 & \dots & 0 & \delta_1 a_N(t)\\ 
0 & - a_2(t) & 0 & \dots & 0 & \delta_2 a_N(t)\\
\cdot & \cdot & \cdot & \dots & \cdot & \cdot \\
0 & 0 & \dots & 0 & - a_{N-1}(t) &  \delta_{N-1} a_N(t)\\
\frac{1}{Z} a_1(t) & \frac{1}{Z} a_2(t) & \frac{1}{Z} a_3(t) & \dots &  \frac{1}{Z} a_{N-1}(t) & - a_N(t)\\      
\end{array}
\right)
\end{equation}
where $a_i \in \mathbb{L}^\infty([0,T], E^{\Phi_2}(\mu)), i=1, \dots, N$,
are all nonnegative. 
Assume the initial datum $\vec{f} \in (\mathbb{L}^2(\mu))^N$ 
is nonnegative.
Then the solution $\vec{u}$ of \ref{eq:matrix_cornerstone} is nonnegative.
Moreover, one has 
$$
\sum_{i=1}^{N-1} u_i(t) + Z u_N(t) 
= e^{tL} \Bigl( \sum_{i=1}^{N-1} f_i + Z f_N \Bigr)
$$
and consequently, 
provided $\vec{f} \in E^{\Phi_{2 \theta}}(\mu)$, then
$
\vec{u} \in \mathbb{L}^\infty([0,T],E^{\Phi_{2 \theta}}(\mu)).
$  
\end{elem}

It is easy to check that $v(t) \equiv \sum_{i=1}^{N-1} u_i(t) + Z u_N(t)$
satisfies $\partial_t v = L v(t)$.
We detail a bit positivity argument (the remaining is similar
to the two-by-two case).

Let $\vec{v}$ be the unique weak solution of problem
\ref{eq:matrix_cornerstone_+} with initial condition $\vec{f}$. 
We now show $\vec{v}$ is nonnegative and so it coincides to the unique
solution of \ref{eq:matrix_cornerstone} with initial condition $\vec{f}$. 
Thanks to Steklov calculus, the following computation is made rigorous. We
focus on the last component (which is the most complicated one). Let 
$v_N^- \equiv \max(-v_N,0)$.
One has
\begin{align*}
\frac{1}{2} \frac{d}{dt} \mu( (v_N^-)^2 ) & = - \mu( v_N^- \, \partial_t v_N)\\
 & = - \mu(v_N^- L v_N) + \mu( a_N(t) v_N^+ v_N^-) - \mu(\sum_{i=1}^{N-1} 
\frac{a_i}{Z} v_i^+ \, v_N^- ).
\end{align*}
Fisrt, $- \mu(v_N^- L v_N)= - \mathcal{E}((-v_N)^+, -v_N) \leq 0$
as for any $\su \in \mathcal{D}$, 
$0 \leq \mathcal{E}(\su_+,\su_+) \leq  \mathcal{E}(\su_+,\su)$.
Secondly, $ \mu( a_N v_N^+ v_N^-)=0$. And the third term is trivially
nonpositive as the $a_i$'s are assumed nonnegative.
Hence, $\mu( (v_N^-)^2 ) \leq \mu( (f_N^-)^2) = 0$.


\bigskip

We can state the following theorem.
\begin{ethm}  
\label{thm:extension}
Let $L_i$, $i=1,\dots q$, be Markov generators satisfying
the {\it linearity assumptions} described before. Assume 
the {\it nonlinearity assumptions} 
are satisfied as well and that $\vec{f} \geq 0$ belongs 
to $E^{\Phi_{2 \theta}}(\mu)$, with $\theta$ as in the 
{\it initial data assumption}.  
 
\smallskip
\noindent Then, for any reaction rates $\lambda_i >0$,  
there exists a unique nonnegative
weak solution $\vec{u}$ of problem {\it \eqref{eq:RDP_extension}} 
on $[0, \infty)$. 
\end{ethm}

\color{black}
\color{black}

\begin{appendix}
\section*{Appendix}

\section{The entropic inequality}
\label{appendix_entropic}
 Let $\mu$ be a probability measure.
Let $f \geq 0$ be a measurable function s.t. $f \neq 0$ $\mu$-a.e.
Then the two following assertions are equivalent:
\begin{enumerate}
 \item[i.] $f \in \mathbb{L}^1(\mu)$ and 
$f \log \left( \frac{f}{\mu(f)} \right) \in \mathbb{L}^1(\mu)$,
\item[ii.] $f \log_+ f  \in \mathbb{L}^1(\mu)$.
\end{enumerate}
  
Let us extend $\mathbb{L}^1(\mu)$ to the space $\mathbb{L}^{1,-}_{\text{ext}}(\mu)$ of measurable functions $f$
such that $\mu(f_+) < +\infty$ and define $\mu(f) \equiv \mu(f_+) - \mu(f_-) \in \mathbb{R} \cup \{-\infty\}$
if $f\in \mathbb{L}^{1,-}_{\text{ext}}(\mu)$. (Define also symmetrically $\mathbb{L}^{1,+}_{\text{ext}}(\mu)$).
Note that  $f\in \mathbb{L}^{1,-}_{\text{ext}}(\mu)$
and $g \in  \mathbb{L}^1(\mu)$ implies $f+g \in \mathbb{L}^{1,-}_{\text{ext}}(\mu)$ and 
$\mu(f+g) = \mu(f) + \mu(g)$. Moreover, for any $f,g \in \mathbb{L}^{1,-}_{\text{ext}}(\mu)$, $f \leq g$ implies $\mu(f) \leq \mu(g)$.

\begin{elem}[Entropic inequality] 
Let $\mu$ be a probability measure and let $f$ and $g$ be two measurable functions.
Assume $f \geq 0$ (excluding $f = 0$ $\mu$-a.e.) such that $f \log_+ f  \in \mathbb{L}^1(\mu)$ and $\mu(e^{\gamma g}) < +\infty$
for some $\gamma >0$.  
Then $fg \in \mathbb{L}^{1,-}_{\text{ext}}(\mu)$ and
\begin{equation}
\label{eq:entropic_inequality}
\mu \left(f g \right) \leq \frac{1}{\gamma} 
\, \mu \left(f \, \log \frac{f}{\mu(f)}\right) + \frac{\mu (f)}{\gamma}
          \log \mu \left(e^{\gamma g} \right)
\end{equation}
in $\mathbb{R} \cup \{ - \infty \}$.
\end{elem}

The proof is based on the following inequality 
$\forall x \in \mathbb{R}_+, \forall y \in \mathbb{R},$ $x \, y \leq x \log x - x + e^y.$


\section{Basics on Orlicz spaces}
\label{appendix:Orlicz}

Classical properties of Orlicz spaces can be found in~\cite{rao-ren}.
\subsubsection*{Young functions}
Let $\Phi$ be a Young function, that is $\Phi : \dR \rightarrow \dR$ convex, even such that
$\Phi(0)= 0$ and $\Phi$ is not constant. Note that from this, it follows that $\Phi(x) \geq 0$, 
that $\Phi(x) \rightarrow + \infty $ when $x \rightarrow \infty$ and that $\Phi$ is an
increasing function on $[0, +\infty)$.


\subsubsection*{Associated Orlicz spaces}
The space $\mathbb{L}^\Phi(\mu) = \{ u \in \mathbb{L}^0(\mu) : \exists \varepsilon>0
\text{ s.t. } \mu (\Phi (\varepsilon u)) < \infty \}$ is a vector subspace of $\mathbb{L}^0(\mu)$.


\subsubsection*{Gauge norm}
Let $B_\Phi = \{ u \in \mathbb{L}^0(\mu) : 
\mu ({\Phi(u)}) \leq 1
 \}$. Then $B_\Phi$  is a symmetric ($B_\Phi = - B_\Phi$) convex set in $\mathbb{L}^\Phi(\mu)$ 
containing $0$ and satisfying 
\begin{equation}
 \label{eq:L_Phi_covering}
\mathbb{L}^\Phi(\mu) = \cup_{\lambda >0} \lambda B_\Phi.
\end{equation}


From these properties, it follows that the gauge norm 
$$
\|u\|_\Phi \equiv \inf \{ \lambda >0 : u \in \lambda B_\Phi \}
$$
associated to $B_\Phi$ is indeed a norm. One has
\begin{equation}
\label{eq:inverse_of_norm}
 \|u\|_\Phi^{-1} = \sup \{ \gamma >0 : \mu( \Phi( \gamma u )) \leq 1 \}.
\end{equation}
The space $(\mathbb{L}^\Phi(\mu),\| \cdot \|_\Phi)$ is a Banach space.

\subsubsection*{Comparison of norms}

We often have to compare Orlicz norms
associated to different Young functions. We already have seen in a footnote that any
Young function $\Phi$ satisfies $ \left| x \right| \preceq \Phi(x)$. 
It leads to the following
lemma. 

\begin{elem}
\label{lemma_injection_L1} Any Orlicz space may be continuously
embedded in $\mathbb{L}_1(\mu)$. More precisely, let $M$ and $\tau$ in
$(0,\infty)$ such that $ \left| x \right| \leq \tau \,  \Phi(x) $
for any $\left| x \right| \geq M$. Then, for any $f \in \mathbb{L}_\Phi$,
\begin{equation}
\label{norm_injection} \NRM{f}_1 \leq (M+\tau) \, \NRM{f}_\Phi.
\end{equation}
Consequently, if $\Phi$ and $\Psi$ are two Young functions
satisfying, for some constants $A,B \geq 0$, $\Phi(x) \leq A  |x| +
B  \Psi(x)$, then
\begin{equation}
\label{norm_comparison} \NRM{f}_\Phi \leq \max \left(1 ,  A
\NRM{\text{{\it Id}}}_{\mathbb{L}_\Psi
    \rightarrow \mathbb{L}_1} + B  \right)  \NRM{f}_\Psi.
\end{equation}
\end{elem}

\begin{edefi}[Comparison of Young functions]
Let us denote $\Phi(x) \preceq \widetilde{\Phi}(x)$ if there exist $x_0 \geq 0$ and $C\in(0,+\infty)$ 
such that $\forall x \geq x_0$, $\Phi(x) \leq C \widetilde{\Phi}(x)$. Furthermore,
$\Phi(x) \simeq \widetilde{\Phi}(x)$ will mean $\Phi(x) \preceq \widetilde{\Phi}(x)$
and $\widetilde{\Phi}(x) \preceq \Phi(x)$.
\end{edefi}

\begin{erem}
Let $\Phi$ and $\widetilde{\Phi}$ be two Young
functions. The existence of a constant $A$ such that
$$
\forall x \geq 0, \Phi(x) \leq A \, \left( |x| + \widetilde{\Phi} (x) \right) 
$$ 
is equivalent to the comparison
$$
\Phi(x) \preceq \widetilde{\Phi}(x). 
$$
The previous lemma then claims briefly that comparison of Young functions
induces comparison of norms.
\end{erem}

Indeed, first assume $\forall x \geq 0$, 
$\Phi(x) \leq A \, \left( |x| + \widetilde{\Phi} (x) \right)$. 
As $|x| = \text{O} (\widetilde{\Phi}(x))$  as  $x$ goes to $+\infty$,
there exist $x_0 \text{ and } B \text{ s.t. } \forall x \geq x_0,
|x| \leq B \widetilde{\Phi} (x)$. So that 
$\forall x \geq x_0, \Phi(x) \leq A (B + 1) \widetilde{\Phi} (x)$.

\smallskip
Conversely, $\Psi(x) \equiv |x| + \widetilde{\Phi}(x)$ is a
Young function, so that $\frac{\Psi(x)}{x}$ is non decreasing
on $(0,\infty)$ and $\forall x > 0, \frac{\Psi(x)}{x} \geq \Psi'(0_+) \geq 1$.
Hence, for any $0 <x \leq x_0$,
$$
\frac{\Phi(x)}{\Psi(x)}=  \frac{\Phi(x)}{x} \underbrace{\frac{x}{\Psi(x)}}_{\leq 1} 
\leq  \frac{\Phi(x)}{x} \leq  \frac{\Phi(x_0)}{x_0}.
$$
The result follows with $A= \max (C, \frac{\Phi(x_0)}{x_0})$.

\medskip
We will also need to deduce bounds on conjugate functions (as defined 
in \eqref{eq:conjugate_function})
from bounds on Young functions. 

\begin{elem}[\cite{rao-ren}, Proposition II.2]
\label{lemma:rao-ren}
Let $\Phi$ and $\Psi$ be  Young functions 
and $\Phi^\ast$ and $\Psi^\ast$ their conjugate functions.
Assume there exits $x_0 \geq 0$ such that
$$
\forall x \geq x_0, \quad \Phi(x) \leq \Psi(x).
$$
Then, there exists $y_0 \geq 0$ such that
$$
\forall y \geq y_0, \quad \Psi^\ast(y) \leq \Phi^\ast(y).
$$ 
\end{elem}

\subsubsection*{Exponential type Young functions and their conjugates}
Let us recall we considered Young functions of exponential type
$$
\Phi_\alpha(x) = \exp(|x|^\alpha) -1, \qquad \alpha \geq 1.
$$
A direct computation shows that, for $y \geq 0$, 
$$
\Phi_1^\ast (y) = 
\left\{ \begin{array}{ll}
         0 & \text{ if } y \leq 1\\
         y \log y - y + 1 &  \text{ if } y \geq 1 
        \end{array}
\right. 
$$
As a consequence, $\Phi_1^\ast (y) \simeq h(y) \equiv y \log_+ y$ 
and $\Phi_1^\ast$ is $\Delta_2$. Here $\log_+ y= \max(\log y, 0)$.
Using lemmas~\ref{lemma:rao-ren} and~\ref{lemma_injection_L1}, 
it follows that, provided $1 \leq \alpha \leq \beta < \infty$
\begin{equation}
\label{eq:comparison-conjugates}
\Phi_\beta^\ast \preceq  \Phi_\alpha^\ast \preceq h \preceq x^2 
\text{ so that } \qquad \| \cdot \|_{\Phi_\beta^\ast} \preceq 
\| \cdot \|_{\Phi_\alpha^\ast} \preceq \| \cdot \|_{h} \preceq \| \cdot \|_{2} 
\end{equation}

\subsubsection*{More on $\mathbf{E^{\Phi_\alpha}(\mu)}$}

One may change parameters in Young inequality to get: 
for any $\alpha >1$ and any $\delta, r >0$, one has
$ 
\forall s \geq 0, \qquad 
\exp( \delta s ) \leq \exp( \frac{\alpha - 1}{\alpha}
( r \alpha / \delta^\alpha )^{\frac{1}{1-\alpha}}) \, \exp( r s^\alpha).
$
It follows that, for any $\alpha \geq 1$, 
$$
\cup_{\beta >\alpha} \mathbb{L}^{\Phi_\beta}(\mu) \subset E^{\Phi_\alpha}(\mu).
$$
 
\begin{elem}[Separability]
Assume $\mathbb{M}$ is a separable metric space. Then, for any Young 
function $\Phi$, $E^\Phi(\mu)$ is separable.
\end{elem}

(Use that $\mathcal{B}_{\mathbb{M}}$ is countably generated,  
monotone class theorem and density of simple functions).

\subsubsection*{Duality}

What follows may be found in \cite{chen-orlicz}.

A Young function $\Psi:\dR \rightarrow \dR$ is said to satisfy the $\Delta_2$ condition
if there exists $K \in (0,\infty)$ and $x_0 \geq 0$ such that, for any $x \geq x_0$,
$\Psi(2x) \leq K \Psi(x)$. 

In the case of Young functions with rapid
growth (as the $\Phi_\alpha$'s introduced before), 
$\Delta_2$ condition fails. Consequently
$E^\Phi(\mu)$ is a proper Banach subspace of $\mathbb{L}^\Phi(\mu)$
(assuming the support of $\mu$ is infinite
) and 
$\mathbb{L}^\Phi(\mu)$ is not separable. 

\medskip
Recall that the conjugate function $\Psi^\ast$ of a Young function $\Psi$ 
is the Young function defined by 
\begin{equation}
\label{eq:conjugate_function}
\Psi^\ast (y) \equiv \sup_{x \geq 0} ( x |y| - \Psi(x)).
\end{equation}
The dual space  of $E^\Psi(\mu)$ is $E^\Psi(\mu)' = \mathbb{L}^{\Psi^\ast}(\mu)$. 

But when 
$\Delta_2$ condition fails, the dual space of $\mathbb{L}^{\Psi}(\mu)$ is more complicated:
this is a direct sum of 
$\mathbb{L}^{\Psi^\ast}(\mu)$ with some nontrivial subspace made of {\it singular} 
linear forms. As a consequence, neither $\mathbb{L}^{\Phi_\alpha}(\mu)$, 
$E^{\Phi_\alpha}(\mu)$ nor $\mathbb{L}^{\Phi_\alpha^\ast}(\mu)$ is reflexive.

\section{Markov Semigroups and Orlicz spaces}

\label{appendix:Orlicz_semigroup}

\subsection{Contraction property}

\begin{elem}
\label{lemma:Jensen_Pt_Orlicz}
Let $\Phi:\dR \rightarrow \dR_+$ be a nonnegative convex function.
Let $(P_t)_{t\geq 0}$ be a Markov semigroup on $\mathbb{L}^2(\mu)$, for
 a probability measure $\mu$,  as introduced
in section~\ref{section:main_result}. Then, for any $f \in \mathbb{L}^1(\mu)$ and
any $t \geq 0$,
\begin{equation}
\label{eq:contraction_of_modulus}
 \mu ( \Phi(P_t f) ) \leq \mu \Phi(f). 
\end{equation}
In particular, in the case when $\Phi$ is a Young function (with domain $\dR$), provided
$f \in \mathbb{L}^\Phi(\mu)$, then
$P_t f \in \mathbb{L}^\Phi(\mu)$ and $(P_t)_{t \geq 0}$ is a
contraction semigroup on $\mathbb{L}^\Phi(\mu)$.
\end{elem}

\begin{eproof}
Let $f \in \mathbb{L}^1(\mu)$, $t \geq 0$ and $\Phi: \dR \rightarrow \dR_+$ be convex.
Nonnegativity of $\Phi$ allows to use Jensen inequality for the Markov
probability kernels $p_t(x,dy)$. Indeed, for $\mu$ almost every $x \in \mathbb{M}$
(such that the representation \eqref{eq:noyau} holds) and any $y \in \mathbb{M}$, by convexity,
$$
\Phi(f(y)) \geq \Phi(P_t f (x))  + \Phi'((P_t f (x))_+) \, (f(y) - P_tf(x)).
$$
Integrating w.r.t. $p_t(x,dy)$ leads to
$$
P_t( \Phi(f) ) (x) \geq \Phi(P_tf(x)) \geq 0. 
$$
Then \eqref{eq:contraction_of_modulus} follows by integration w.r.t. $\mu$
and invariance property of $P_t$.

\bigskip
Let now $\Phi$ be a Young function. Assume $f \neq 0$ in 
$\mathbb{L}^\Phi(\mu) (\subset \mathbb{L}^1(\mu))$.  Recall \eqref{eq:inverse_of_norm} and choose
$0 < \gamma \leq \|f \|_{\mathbb{L}^\Phi}^{-1}$.
Applying \eqref{eq:contraction_of_modulus} to 
$\Phi( \gamma \cdot )$ instead of $\Phi$ shows that
$$
\mu(\Phi(\gamma P_t f)) \leq 1
$$
so that $\gamma \leq  \|P_t f \|_{\mathbb{L}^\Phi}^{-1}$.
And the announced contraction property follows.
\end{eproof}

\subsection{Density of the Dirichlet domain}

Using comparison \eqref{eq:comparison-conjugates}, one gets  
continuous embedding 
\begin{equation}
\label{eq:embedding_Dirichlet_Phiast}
 \mathcal{D} \hookrightarrow \mathbb{L}^2(\mu) \hookrightarrow \mathbb{L}^{\Phi_\alpha^\ast}(\mu),
\end{equation}
for any $\alpha \geq 1$. As $\Phi_\alpha^\ast$ is $\Delta_2$, the space of simple functions, and so
$\mathbb{L}^2(\mu)$ as well, is dense in  $\mathbb{L}^{\Phi_\alpha^\ast}(\mu)$.
Now, $\mathcal{D}$ is dense in $\mathbb{L}^2(\mu)$, and so in  
$\mathbb{L}^{\Phi_\alpha^\ast}(\mu)$. 

\subsection{Bochner measurability}
\label{appendix:Bochner-measurability}

Let $X$ be a Banach space.
Recall that an $X$-valued function $u:I \rightarrow X$ defined on
a compact interval $I$ is Bochner measurable provided it is an $a.e.$
limit of a sequence of $X$-valued simple functions on $I$
(see \cite{schwabik-ye} for instance).

\subsubsection*{The $\mathbf{\mathbb{L}^\infty([0,T],\mathbb{L}^\Phi(\mu))}$ space}

\begin{elem}
\label{lemma:LinftyLPhi}
Let $\Phi:\dR \rightarrow \dR$ be a Young function, $(\mathbb{M},\mathcal{B}_\mathbb{M},\mu)$
a probability space and $u \in C([0,T],\mathbb{L}^2(\mu))$. We assume that 
$x^2 \preceq \Phi(x)$. Then
$
u \in \mathbb{L}^\infty([0,T],\mathbb{L}^\Phi(\mu)) 
$
iff
$ 
u : t \rightarrow u(t,\cdot) \in \mathbb{L}^\Phi(\mu)
$  
is Bochner measurable and  
there exist $\gamma, M \in (0,\infty)$ 
s.t., for any $t$ a.e. in $[0,T]$,
$ 
\mu( \Phi(\gamma u(t))) \leq M.
$
In which case, one has, for any $t$ a.e., 
$\|u(t)\|_{\mathbb{L}^\Phi(\mu)} \leq \frac{\max(M,1)}{\gamma}$.
\end{elem}

This is just rewritting the definitions. In particular,
provided $M \geq 1$,  and 
$ \mu( \Phi(\gamma u(t))) \leq M$, then by convexity,
$ \mu( \Phi(\frac{\gamma}{M} u(t))) 
\leq \frac{1}{M} \mu( \Phi({\gamma} u(t))) \leq 1
$
so that  $\|u(t)\|_{\mathbb{L}^\Phi(\mu)} \leq \frac{M}{\gamma}$.

\subsubsection*{Proof of proposition \ref{proposition:Bochner_measurability_semigroup}}
By density of $\mathbb{L}^2(\mu)$ in $\mathbb{L}^{\Phi_\alpha^\ast}$
and contraction of $P_t$ in $\mathbb{L}^{\Phi_\alpha^\ast}$, $C_0$ property
of $P_t$ in $\mathbb{L}^{\Phi_\alpha^\ast}$ follows from $C_0$ property
in $\mathbb{L}^2(\mu)$. Indeed, let $f \in \mathbb{L}^{\Phi_\alpha^\ast}$.
$\varepsilon >0$ being fixed, let $g \in \mathbb{L}^2(\mu)$
such that $\|f-g\|_{\Phi_\alpha^\ast} < \frac{\varepsilon}{3}$.
Then
$$
\|P_tf -f\|_{\Phi_\alpha^\ast} \leq 
2 \|f-g\|_{\Phi_\alpha^\ast} + \|P_tg -g\|_{\Phi_\alpha^\ast}
\leq \frac{2 \varepsilon}{3} + C \|P_tg -g\|_{2}
$$  
allows to conclude.
As a consequence, provided $f \in E^{\Phi_\alpha}$, 
$t \mapsto P_tf \in E^{\Phi_\alpha}$ is weakly continuous, 
and so Bochner measurable as $E^{\Phi_\alpha}$ is separable, following Pettis
measurability theorem (see page \pageref{pettis} for references).

\end{appendix}

\bigskip
\noindent
{\bf Acknowledgements.} 
This research was supported  by the ANR project STAB; by the GDR project
AFHP  and the ESP team of IMT; B.Z. was supported 
by Royal Society Wolfson RMA.

\bibliographystyle{plain}

\begin{thebibliography}{}

\end{thebibliography}


\begin{thebibliography}{10}

\bibitem{amann-local}
H.~Amann.
\newblock Existence and regularity for semilinear parabolic evolution
  equations.
\newblock {\em Ann. Scuola Norm. Sup. Pisa Cl. Sci. (4)}, 11(4):593--676, 1984.

\bibitem{amann-global}
H.~Amann.
\newblock Global existence for semilinear parabolic systems.
\newblock {\em J. Reine Angew. Math.}, 360:47--83, 1985.

\bibitem{ancona}
A.~Ancona.
\newblock Continuit\'e des contractions dans les espaces de {D}irichlet.
\newblock In {\em S\'eminaire de {T}h\'eorie du {P}otentiel de {P}aris, {N}o. 2
  ({U}niv. {P}aris, {P}aris, 1975--1976)}, pages 1--26. Lecture Notes in Math.,
  Vol. 563. Springer, Berlin, 1976.

\bibitem{bakry-emery}
D.~Bakry and Michel {\'E}mery.
\newblock Diffusions hypercontractives.
\newblock In {\em S\'eminaire de probabilit\'es, {XIX}, 1983/84}, volume 1123
  of {\em Lecture Notes in Math.}, pages 177--206. Springer, Berlin, 1985.

\bibitem{bakry-gentil-ledoux}
D.~Bakry, I.~Gentil, and M.~Ledoux.
\newblock {\em Analysis and Geometry of Markov Diffusion Operators}, volume 348
  of {\em Grundlehren der mathematischen Wissenschaften}.
\newblock Springer, 2014.

\bibitem{barthe-cattiaux-roberto}
F.~Barthe, P.~Cattiaux, and C.~Roberto.
\newblock Interpolated inequalities between exponential and {G}aussian,
  {O}rlicz hypercontractivity and isoperimetry.
\newblock {\em Rev. Mat. Iberoam.}, 22(3):993--1067, 2006.

\bibitem{bobkov-zegarlinski-slow-tails}
S.~Bobkov and B.~Zegarlinski.
\newblock Distributions with slow tails and ergodicity of {M}arkov semigroups
  in infinite dimensions.
\newblock In {\em Around the research of {V}ladimir {M}az'ya. {I}}, volume~11
  of {\em Int. Math. Ser. (N. Y.)}, pages 13--79. Springer, New York, 2010.

\bibitem{bobkov-zegarlinski}
S.~G. Bobkov and B.~Zegarlinski.
\newblock Entropy bounds and isoperimetry.
\newblock {\em Mem. Amer. Math. Soc.}, 176(829):x+69, 2005.

\bibitem{bodineau-helffer}
T.~Bodineau and B.~Helffer.
\newblock The log-{S}obolev inequality for unbounded spin systems.
\newblock {\em J. Funct. Anal.}, 166(1):168--178, 1999.

\bibitem{bouleau-hirsch}
N.~Bouleau and F.~Hirsch.
\newblock {\em Dirichlet forms and analysis on {W}iener space}, volume~14 of
  {\em de Gruyter Studies in Mathematics}.
\newblock Walter de Gruyter \& Co., Berlin, 1991.

\bibitem{carrillo-hittmeir-jungel}
J.~A. Carrillo, S.~Hittmeir, and A.~J{\"u}ngel.
\newblock Cross diffusion and nonlinear diffusion preventing blow up in the
  {K}eller-{S}egel model.
\newblock {\em Math. Models Methods Appl. Sci.}, 22(12):1250041, 35, 2012.

\bibitem{chen}
M.-F. Chen.
\newblock {\em Eigenvalues, inequalities, and ergodic theory}.
\newblock Probability and its Applications (New York). Springer-Verlag London,
  Ltd., London, 2005.

\bibitem{chen-orlicz}
S.~Chen.
\newblock Geometry of {O}rlicz spaces.
\newblock {\em Dissertationes Math. (Rozprawy Mat.)}, 356:204, 1996.
\newblock With a preface by Julian Musielak.

\bibitem{davies90}
E.~B. Davies.
\newblock {\em Heat kernels and spectral theory}, volume~92 of {\em Cambridge
  Tracts in Mathematics}.
\newblock Cambridge University Press, Cambridge, 1990.

\bibitem{desvillettes-seul}
L.~Desvillettes.
\newblock About entropy methods for reaction-diffusion equations.
\newblock {\em Riv. Mat. Univ. Parma (7)}, 7:81--123, 2007.

\bibitem{df06}
L.~Desvillettes and K.~Fellner.
\newblock Exponential decay toward equilibrium via entropy methods for
  reaction-diffusion equations.
\newblock {\em J. Math. Anal. Appl.}, 319(1):157--176, 2006.

\bibitem{df08}
L.~Desvillettes and K.~Fellner.
\newblock Entropy methods for reaction-diffusion equations: Slowly growing
  {A}-priori bounds.
\newblock {\em Rev. Mat. Iberoamericana}, 24(2):407--431, 2008.

\bibitem{diaconis-saloff}
P.~Diaconis and L.~Saloff-Coste.
\newblock Logarithmic {S}obolev inequalities for finite {M}arkov chains.
\newblock {\em Ann. Appl. Probab.}, 6(3):695--750, 1996.

\bibitem{diestel-uhl}
J.~Diestel and J.~J. Uhl, Jr.
\newblock {\em Vector measures}.
\newblock American Mathematical Society, Providence, R.I., 1977.
\newblock With a foreword by B. J. Pettis, Mathematical Surveys, No. 15.

\bibitem{zelik04}
M.~Efendiev, A.~Miranville, and S.~V. Zelik.
\newblock Infinite-dimensional exponential attractors for nonlinear
  reaction-diffusion systems in unbounded domains and their approximation.
\newblock {\em Proc. R. Soc. Lond. Ser. A Math. Phys. Eng. Sci.},
  460(2044):1107--1129, 2004.

\bibitem{evans}
L.~C. Evans.
\newblock {\em Partial differential equations}, volume~19 of {\em Graduate
  Studies in Mathematics}.
\newblock American Mathematical Society, Providence, RI, second edition, 2010.

\bibitem{FRZ}
P.~Foug{\`e}res, C.~Roberto, and B.~Zegarlinski.
\newblock Sub-gaussian measures and associated semilinear problems.
\newblock {\em Rev. Mat. Iberoam.}, 28(2):305--350, 2012.

\bibitem{friedman}
A.~Friedman.
\newblock {\em Partial Differential Equations of Parabolic Type}.
\newblock Prentice Hall, Englewood Cliffs, N.J., 1964.

\bibitem{fukushima-al}
M.~Fukushima, Y.~Oshima, and M.~Takeda.
\newblock {\em Dirichlet forms and symmetric {M}arkov processes}, volume~19 of
  {\em de Gruyter Studies in Mathematics}.
\newblock Walter de Gruyter \& Co., Berlin, extended edition, 2011.

\bibitem{gentil-zeg}
I.~Gentil and B.~Zegarlinski.
\newblock Asymptotic behaviour of reversible chemical reaction-diffusion
  equations.
\newblock {\em Kinet. Relat. Models}, 3(3):427--444, 2010.

\bibitem{gross75}
L.~Gross.
\newblock Logarithmic {S}obolev inequalities.
\newblock {\em Amer. J. Math.}, 97(4):1061--1083, 1975.

\bibitem{guionnet-zegarlinski}
A.~Guionnet and B.~Zegarlinski.
\newblock Lectures on logarithmic {S}obolev inequalities.
\newblock In {\em S\'eminaire de {P}robabilit\'es, {XXXVI}}, volume 1801 of
  {\em Lecture Notes in Math.}, pages 1--134. Springer, Berlin, 2003.

\bibitem{hebisch-zegarlinski}
W.~Hebisch and B.~Zegarli{\'n}ski.
\newblock Coercive inequalities on metric measure spaces.
\newblock {\em J. Funct. Anal.}, 258(3):814--851, 2010.

\bibitem{inglis-papageorgiou}
J.~Inglis and I.~Papageorgiou.
\newblock Logarithmic {S}obolev inequalities for infinite dimensional
  {H}\"ormander type generators on the {H}eisenberg group.
\newblock {\em Potential Anal.}, 31(1):79--102, 2009.

\bibitem{ladyzenskaja}
O.~A. Lady{\v{z}}enskaja, V.~A. Solonnikov, and N.~N. Ural{\cprime}ceva.
\newblock {\em Linear and quasilinear equations of parabolic type}.
\newblock Translated from the Russian by S. Smith. Translations of Mathematical
  Monographs, Vol. 23. American Mathematical Society, Providence, R.I., 1968.

\bibitem{ledoux-spins}
M.~Ledoux.
\newblock Logarithmic {S}obolev inequalities for unbounded spin systems
  revisited.
\newblock In {\em S\'eminaire de {P}robabilit\'es, {XXXV}}, volume 1755 of {\em
  Lecture Notes in Math.}, pages 167--194. Springer, Berlin, 2001.

\bibitem{lugiewicz-zegarlinski}
P.~{\L}ugiewicz and B.~Zegarli{\'n}ski.
\newblock Coercive inequalities for {H}\"ormander type generators in infinite
  dimensions.
\newblock {\em J. Funct. Anal.}, 247(2):438--476, 2007.

\bibitem{ma-rockner}
Z.~M. Ma and M.~R{\"o}ckner.
\newblock {\em Introduction to the theory of (nonsymmetric) {D}irichlet forms}.
\newblock Universitext. Springer-Verlag, Berlin, 1992.

\bibitem{mahe-fraissard}
R.~Mah\'e and J.~Fraissard.
\newblock {\em \'Equilibres chimiques en solution acqueuse}.
\newblock Masson, Paris, 1989.

\bibitem{miclo}
L.~Miclo.
\newblock An example of application of discrete {H}ardy's inequalities.
\newblock {\em Markov Process. Related Fields}, 5(3):319--330, 1999.

\bibitem{otto-reznikoff}
F.~Otto and M.~G. Reznikoff.
\newblock A new criterion for the logarithmic {S}obolev inequality and two
  applications.
\newblock {\em J. Funct. Anal.}, 243(1):121--157, 2007.

\bibitem{pierre}
M.~Pierre.
\newblock Global existence in reaction-diffusion systems with control of mass:
  a survey.
\newblock {\em Milan J. Math.}, 78(2):417--455, 2010.

\bibitem{rao-ren}
M.~M. Rao and Z.~D. Ren.
\newblock {\em Theory of {O}rlicz spaces}, volume 146 of {\em Monographs and
  Textbooks in Pure and Applied Mathematics}.
\newblock Marcel Dekker Inc., New York, 1991.

\bibitem{roberto-zegarlinski}
C.~Roberto and B.~Zegarli{\'n}ski.
\newblock Orlicz-{S}obolev inequalities for sub-{G}aussian measures and
  ergodicity of {M}arkov semi-groups.
\newblock {\em J. Funct. Anal.}, 243(1):28--66, 2007.

\bibitem{rothe}
F.~Rothe.
\newblock {\em Global solutions of reaction-diffusion systems}, volume 1072 of
  {\em Lecture Notes in Mathematics}.
\newblock Springer-Verlag, Berlin, 1984.

\bibitem{schwabik-ye}
S.~Schwabik and G.~Ye.
\newblock {\em Topics in {B}anach space integration}, volume~10 of {\em Series
  in Real Analysis}.
\newblock World Scientific Publishing Co. Pte. Ltd., Hackensack, NJ, 2005.

\bibitem{stroock-zegarlinski}
D.~W. Stroock and B.~Zegarli{\'n}ski.
\newblock The logarithmic {S}obolev inequality for continuous spin systems on a
  lattice.
\newblock {\em J. Funct. Anal.}, 104(2):299--326, 1992.

\bibitem{taylorIII}
M.~E. Taylor.
\newblock {\em Partial differential equations. {III}}, volume 117 of {\em
  Applied Mathematical Sciences}.
\newblock Springer-Verlag, New York, 1997.
\newblock Nonlinear equations, Corrected reprint of the 1996 original.

\bibitem{wang}
F.Y. Wang.
\newblock {\em Functional Inequalities, Markov Semigroups and Spectral Theory}.
\newblock Science Press, Beijing, 2005.

\bibitem{yoshida00}
N.~Yoshida.
\newblock Application of log-{S}obolev inequality to the stochastic dynamics of
  unbounded spin systems on the lattice.
\newblock {\em J. Funct. Anal.}, 173(1):74--102, 2000.

\bibitem{yosida}
K.~Yosida.
\newblock {\em Functional analysis}.
\newblock Die Grundlehren der Mathematischen Wissenschaften, Band 123. Academic
  Press, Inc., New York; Springer-Verlag, Berlin, 1965.

\bibitem{zegarlinski90}
B.~Zegarli{\'n}ski.
\newblock On log-{S}obolev inequalities for infinite lattice systems.
\newblock {\em Lett. Math. Phys.}, 20(3):173--182, 1990.

\bibitem{zegarlinski96}
B.~Zegarlinski.
\newblock The strong decay to equilibrium for the stochastic dynamics of
  unbounded spin systems on a lattice.
\newblock {\em Comm. Math. Phys.}, 175(2):401--432, 1996.

\bibitem{zelik03}
S.~V. Zelik.
\newblock Attractors of reaction-diffusion systems in unbounded domains and
  their spatial complexity.
\newblock {\em Comm. Pure Appl. Math.}, 56(5):584--637, 2003.

\bibitem{zelik07}
S.~V. Zelik.
\newblock Spatial and dynamical chaos generated by reaction-diffusion systems
  in unbounded domains.
\newblock {\em J. Dynam. Differential Equations}, 19(1):1--74, 2007.

\end{thebibliography}

\def\cprime{$'$}

\noindent
Institut de Math\'ematiques de Toulouse, UMR CNRS 5219\\
Universit\'e de Toulouse\\
Route de Narbonne\\
31062 Toulouse - France\\
fougeres@math.univ-toulouse.fr

\medskip

\noindent
Institut Camille Jordan, UMR CNRS 5208\\
Universit\'e Claude Bernard, Lyon 1\\
43 boulevard du 11 novembre 1918\\
69622 Villeurbanne cedex - France\\
gentil@math.univ-lyon1.fr

\medskip

\noindent
Imperial College, London\\
South Kensington Campus\\
London SW7 2AZ\\
b.zegarlinski@imperial.ac.uk

\end{document}